\documentclass[
               DIV=15,
               ]{scrartcl}  

\usepackage{color} 
\usepackage{psfrag} 
\usepackage{subcaption}
\usepackage{hyperref}
\usepackage{amsmath}
\usepackage{amssymb}
\usepackage{amsthm}
\usepackage{graphics}
\usepackage{todonotes}
\usepackage{amsfonts}  
\usepackage{adjustbox}    
\usepackage{mathtools}
\usepackage{xcolor}        
\usepackage{blkarray}        
\usepackage{footnote}
\makesavenoteenv{tabular}
\makesavenoteenv{table}

\usepackage[square,numbers,comma]{natbib}

\hypersetup{colorlinks,citecolor=blue}

\usepackage{standalone}
 
\usepackage{pgfplots}
\usepackage{pgfplotstable}
\usetikzlibrary{pgfplots.groupplots}
\usetikzlibrary{fit}
\usepackage{filecontents}
\usepackage{tikz-3dplot}

\usepgfplotslibrary{patchplots}

\usepackage{tikz}
\usepackage{tikzscale}

\usetikzlibrary{arrows,3d}

\usetikzlibrary{shapes.geometric}
\usetikzlibrary{positioning}
\usetikzlibrary{calc}
\usetikzlibrary{decorations.pathreplacing} 
\usetikzlibrary{decorations.markings}

\usetikzlibrary{spy}
\usetikzlibrary{snakes,arrows,shapes}    
\usetikzlibrary{spy}
\usetikzlibrary{backgrounds,tikzmark}  
\usetikzlibrary{decorations}
\usetikzlibrary{positioning}
\usetikzlibrary{patterns}
\usetikzlibrary{calc} 

\usepackage{booktabs}   
\usepackage{makecell} 
\usepackage{colortbl}   

\usepackage{xspace}
\usepackage{color}     
\usepackage{setspace}   
 
\usepackage{soul}
\usepackage{ulem}
\usepackage{currfile}
\usepackage{import}

\usepackage{authoraftertitle} 
\usepackage{authblk} 

\usepackage[noabbrev]{cleveref}

\usepackage{enumitem}
\usepackage{multicol} 
\usepackage{algorithm}
\usepackage{algpseudocode}
\usepackage{algorithmicx}
\usepackage[section]{placeins}

\usepackage{environ}

    \usepackage{framed}
    
    \setlist{noitemsep,topsep=2pt,parsep=0pt,partopsep=0pt}
  
    \usepackage{lmodern}

    \usetikzlibrary{external}
    \tikzset{external/force remake}
    
   \pgfdeclarelayer{background}
   \pgfdeclarelayer{foreground}
   \pgfsetlayers{background,main,foreground}

\pgfplotsset{compat=1.8}

\newcommand{\findmax}[3]{
    \pgfplotstablesort[sort key={#2},sort cmp={float >}]{\sorted}{#1}%
    \pgfplotstablegetelem{0}{#2}\of{\sorted}%
    \let #3=\pgfplotsretval%
}

\pgfplotscreateplotcyclelist{myCycleListBlackAndWhite}
{%
  {black, mark=o},
  {black, mark=square},
  {black, mark=triangle},
  {black, mark=diamond}, 
  {black, mark=pentagon}, 
  {black, mark=*},
  {black, mark=square*},
  {black, mark=triangle*},
  {black, mark=diamond*}, 
  {black, mark=pentagon*}, 
}

\definecolor{darkgreen}{rgb}{0,0.4,0} 
\definecolor{darkbrown}{rgb}{0.5, 0.396, 0.09}

\definecolor{c1}{rgb}{0.0, 0.4196078431372549, 0.6431372549019608}
\definecolor{c2}{rgb}{1.0, 0.5019607843137255, 0.054901960784313725}
\definecolor{c3}{rgb}{0.6705882352941176, 0.6705882352941176,
0.6705882352941176} \definecolor{c}{rgb}{0.34901960784313724, 0.34901960784313724, 0.34901960784313724}
\definecolor{c4}{rgb}{0.37254901960784315, 0.6196078431372549,
0.8196078431372549} \definecolor{c}{rgb}{0.7843137254901961, 0.3215686274509804, 0.0}
\definecolor{c5}{rgb}{0.5372549019607843, 0.5372549019607843,
0.5372549019607843} \definecolor{c}{rgb}{0.6352941176470588, 0.7843137254901961, 0.9254901960784314}
\definecolor{c6}{rgb}{1.0, 0.7372549019607844, 0.4745098039215686}
\definecolor{c7}{rgb}{0.8117647058823529, 0.8117647058823529,
0.8117647058823529}

\pgfplotscreateplotcyclelist{myCycleListColor}
{%
  {blue, mark=o},
  {red, mark=square},
  {darkbrown, mark=triangle},
  {darkgreen, mark=diamond},
  {black, mark=pentagon},
  {blue, mark=*},
  {red, mark=square*},
  {darkbrown, mark=triangle*},
  {darkgreen, mark=diamond*},
  {black, mark=pentagon*},
}

\pgfplotsset{every axis/.append style= 
              {
                font=\scriptsize,
                mark size=3,
                legend style={font=\tiny, mark size=3, draw=none, fill=none},
                legend cell align=left,
                cycle list name=myCycleListColor,
                scaled y ticks = false,
				scaled x ticks = false,
				trim axis left,
				trim axis right,
				sharp plot,
				tick label style ={font=\tiny},
				label style ={font=\scriptsize},
				very thin,
				ymajorgrids=true,
				grid style=dotted,
				legend pos= north east,
				height=\figHeight, width=\figWidth,
              }
            }

\pgfplotsset{every tick label/.append style={font=\tiny}}
\pgfplotsset{every label/.append style={font=\scriptsize}}
\pgfplotsset{every axis legend/.append style={font=\tiny}}
\pgfplotsset{every axis plot/.append style={thick}}

\pgfplotstableset
{
  every odd row/.style={before row={\rowcolor[gray]{0.9}}},
  every head row/.style=
  {
    before row=
    {%
      \toprule
    },
    after row=\midrule
  },
  every last row/.style=
  {
    after row=\bottomrule
  }
}

\newif\ifdrawboundingbox

\pgfkeys
{ 
  /pgf/number format/.cd,
  fixed,
  set thousands separator={\,},
  1000 sep in fractionals,
}

\newcommand{\picsDir}{pictures/numericalExamples/pics}

\newcolumntype{C}[1]{>{\centering\arraybackslash}m{#1}}
\newcolumntype{R}[1]{>{\raggedright\arraybackslash}m{#1}}
\newcolumntype{L}[1]{>{\raggedleft\arraybackslash}m{#1}}

\newcommand\restr[2]{{
		\left.\kern-\nulldelimiterspace 
		#1 
		\right|_{#2} 
}}

\newcommand\compactDots{\ifmmode\ldots\else\makebox[10cm][c]{.\hfil.\hfil.}\fi}

\newcommand{\supp}[1]{\textnormal{supp}(#1)}

\def\registered{{\ooalign{\hfil\raise .00ex\hbox{\tiny R}\hfil\cr\mathhexbox20D}}}

\makeatletter
\def\namedlabel#1#2{\begingroup
	\def\@currentlabel{#2}%
	\label{#1}\endgroup
}
\makeatother

\makeatletter
\renewcommand{\todo}[2][]{\tikzexternaldisable\@todo[#1]{#2}\tikzexternalenable}
\makeatother

\usepackage{etoolbox}
\AtBeginEnvironment{mysmallmatrix}{%
\setstretch{0.9}%
\setlength{\arraycolsep}{2pt}%
}%

{%
	\begin{matrix}%
}%
	{\end{matrix}}%

{%
	\begin{matrix}%
	}%
	{\end{matrix}}%
\AtBeginEnvironment{myverysmallmatrix}{%
	\footnotesize%
	\setstretch{0.55}%
	\setlength{\arraycolsep}{1.5pt}%
}%

\makeatletter
\def\bbordermatrix#1{\begingroup \m@th
	\@tempdima 4.75\p@
	\setbox\z@\vbox{%
		\def\cr{\crcr\noalign{\kern2\p@\global\let\cr\endline}}%
		\ialign{$##$\hfil\kern2\p@\kern\@tempdima&\thinspace\hfil$##$\hfil
			&&\quad\hfil$##$\hfil\crcr
			\omit\strut\hfil\crcr\noalign{\kern-\baselineskip}%
			#1\crcr\omit\strut\cr}}%
	\setbox\tw@\vbox{\unvcopy\z@\global\setbox\@ne\lastbox}%
	\setbox\tw@\hbox{\unhbox\@ne\unskip\global\setbox\@ne\lastbox}%
	\setbox\tw@\hbox{$\kern\wd\@ne\kern-\@tempdima\left[\kern-\wd\@ne
		\global\setbox\@ne\vbox{\box\@ne\kern2\p@}%
		\vcenter{\kern-\ht\@ne\unvbox\z@\kern-\baselineskip}\,\right]$}%
	\null\;\vbox{\kern\ht\@ne\box\tw@}\endgroup}
\makeatother

\title{Adaptive isogeometric analysis on \\two-dimensional trimmed domains based on \\a hierarchical approach.}

\author[1]{Luca Coradello\thanks{luca.coradello@epfl.ch, Corresponding Author}}
\author[1]{Pablo Antolin}
\author[1,2]{Rafael V\'{a}zquez}
\author[1,2]{Annalisa Buffa}

\affil[1]{Institute of Mathematics,
 		  \'Ecole Polytechnique F\'ed\'erale de Lausanne, Lausanne, Switzerland.}

\affil[2]{Istituto di Matematica Applicata e Tecnologie Informatiche `E. Magenes' (CNR), Pavia, Italy.}

\newcommand{\publicationDate}{\today}

\usepackage{fancyhdr}
\date{}
\fancypagestyle{plain}{%
  \fancyhf{}%
  \fancyfoot[L]{
  \vspace{-1.25cm} 
  \footnotesize{Preprint submitted to  \hfill
  \publicationDate
  \\
  }} }
 
\begin{document}  
\normalem
\maketitle  
  
\vspace{-1.5cm} 
\hrule 
\section*{Abstract}
The focus of this work is on the development of an error-driven isogeometric framework, capable of automatically performing an adaptive simulation in the context of second- and fourth-order, elliptic partial differential equations defined on two-dimensional \textit{trimmed} domains. The method is steered by an a posteriori error estimator, which is computed with the aid of an auxiliary residual-like problem formulated onto a space spanned by splines with single element support. The local refinement of the basis is achieved thanks to the use of truncated hierarchical B-splines.  
We prove numerically the applicability of the proposed estimator to various engineering-relevant problems, namely the Poisson problem, linear elasticity and Kirchhoff-Love shells, formulated on trimmed geometries. In particular, we study several benchmark problems which exhibit both smooth and singular solutions, where we recover optimal asymptotic rates of convergence for the error measured in the energy norm and we observe a substantial increase in accuracy per-degree-of-freedom compared to uniform refinement. Lastly, we show the applicability of our framework to the adaptive shell analysis of an industrial-like trimmed geometry modeled in the commercial software \textit{Rhinoceros}, which represents the B-pillar of a car.  
\vspace{1cm}

\textit{Keywords}: isogeometric analysis, trimming, truncated hierarchical B-splines, a posteriori error estimator, adaptivity, Kirchhoff-Love shell
\vspace{0.2cm} 
\hrule

\def\Estconst{3}



\section{Introduction}\label{sec:introduction}

Isogeometric Analysis (IGA) was introduced in 2005 in the seminal work~\cite{Hughes2005} and since then has been a successful area of research in the computational mechanics community. By using the same building blocks employed in standard Computer-Aided Design (CAD), namely B-splines, Non-Uniform Rational B-splines (NURBS) and variances thereof, the final goal of IGA is to provide an end-to-end methodology that unifies geometrical design with the analysis. While this is still an open issue, in the last decade an extensive amount of research has been dedicated to IGA in various different fields, we refer to the special issue~\citep{Hughes2017special} for a review of the most prominent work published in recent years. In particular, B-splines based formulations are now built on solid mathematical foundations, e.g. see~\citep{Buffa2014,Bazilevs2006}, and have demonstrated their capabilities in many different areas of engineering, we refer the reader to~\cite{Hughes2005,Cottrell2009}.  

However, to achieve a smooth Design-Through-Analysis workflow at an industrial level, it is imperative that the proper treatment of trimmed models is addressed. We refer to~\citep{Marussig2018} for a review of the current challenges to be faced due to trimming. In this work, we restrict ourselves to surfaces and we tackle the issue at integration level by properly re-parametrizing those elements that are cut, by exploiting the tool presented in~\citep{Antolin2019} for volumetric geometries. This approach is similar to what was introduced in the pioneering work~\citep{Breitenberger2015} and later in~\citep{Guo2018} in the scope of isogeometric shell analysis and in~\citep{Kudela2015} in the context of immersed methods, namely the Finite Cell Method. 

Another important aspect to be highlighted concerns local refinement. It is well-known that the tensor product structure of B-splines hinders the capability of efficiently capture localized features of the solution in small areas of interest. A considerable amount of recent work has been dedicated to overcome this drawback; we mention as potential remedies hierarchical B-splines (HB) \cite{Forsey1988, Greiner1997, Kraft1997} and their variant denoted by truncated hierarchical B-splines (THB) \cite{Giannelli2012, Giannelli2016}, T-splines \cite{Bazilevs2010,Scott2012,Beirao2013} and LR-splines \cite{Dokken2013,Bressan2013}. For efficiency, also the proper coarsening of spline functions constitutes a rich area of research, where we mention as reference~\citep{Lorenzo2017,Garau2018,Hennig2018,Carraturo2019}.

Indeed, the contribution of this paper falls into this realm. The use of locally-refinable splines on complex geometries defined by trimming operations is still in its infancy and, to the best of the authors' knowledge, the definition of suitable error indicators that account for trimming and thus the use of adaptivity in this context has not been subject of thorough research.    

The main goal of this manuscript is the creation of a unified framework for the adaptive isogeometric analysis of PDEs defined on trimmed domains. This task is achieved by extending the authors' previous work on error estimation for linear fourth-order elliptic partial differential equations (PDEs)~\citep{Antolin2019b} and combining it with the trimming tool proposed in~\citep{Antolin2019}. In particular, this constitutes a versatile approach that, potentially, fully automates the analysis of complex trimmed geometries up to a pre-defined error threshold, without the need of any other input from the user. 
We remark that in order to achieve local refinement, we make use of truncated hierarchical B-splines~\citep{Giannelli2012}, where we show that this choice yields a linearly independent basis suitable for the analysis also in the presence of trimming. 
A similar framework was first studied in~\citep{Hoellig_book,Hollig2015,Hollig2012} for hierarchical B-splines defined on a rectangular planar grid, where no surface parametrization is used and local refinement is performed using \textit{a priori} knowledge of the features of the problem. More recently, in~\citep{Marussig2018,Deprenter2019}, the use of truncated hierarchical B-splines in the context of the Laplace problem and linear elasticity was studied, where the emphasis is put on the issues of stability of the basis and bad conditioning of the system matrix caused by trimming and no real error indicator is employed. Lastly in~\citep{Qarariyah2019}, a fourth-order PDE, namely the Biharmonic problem, formulated on implicit domains was analyzed, where again no error-driven adaptive simulation is performed but refinement is steered \textit{a priori} towards geometric features of interest. Our contribution differs from the aforementioned works as our main focus lies on developing a fully adaptive, error-driven numerical framework for second- and fourth-order PDEs defined on trimmed surfaces.
In particular, at first, we validate the proposed method and show its potential by studying the simple Poisson model problem and later we extend it to the analysis of more challenging and engineering-related problems, such as linear elasticity and Kirchhoff-Love shells. Inspired by the work on multi-level estimators presented in~\citep{Bank1993,Vuong2011}, the adaptive procedure is based on the solution of a residual-like variational problem, formulated on a space of splines with compact support, similarly to what was studied by the authors in~\citep{Antolin2019b}. 


Finally, we remark that our method, which handles single-patch geometries for the time being, shows the potential to become a powerful analysis tool in computational mechanics. Indeed, the ability to properly treat trimmed surfaces coming from a commercial CAD software combined with an adaptive loop driven by the proposed error estimator is a substantial step towards a fully automated, end-to-end numerical simulation of complex geometries of industrial relevance. To demonstrate this, we perform an adaptive shell analysis on the B-pillar of a car.  

We structure the paper as follows. \Cref{sec:trimming} introduces the basic framework of isogeometric methods defined on trimmed domains together with a review of (truncated) hierarchical B-splines. Then, \Cref{sec:bubbles} presents the proposed error estimator in the context of trimming and a possible implementation is given. \Cref{sec:numerical_examples} presents several numerical experiments. Finally, we draw some conclusions in \Cref{sec:conclusions}.

\newtheorem{remark}{Remark}

\section{Mathematical framework of trimming} \label{sec:trimming}
\newcommand{\bg}{}


Trimming is a basic mathematical Boolean operation, which allows for an easy description of complex geometries and it is regarded as standard in most commercial CAD softwares. A trimmed surface is composed of two parts, an underlying geometry, that describes the geometric shape, and a set of properly ordered trimming curves, that delimit the regions that are clipped. Despite its simple definition, trimming yields to severe challenges in the context of a smooth Design-Through-Analysis workflow. Indeed, when a trimming operation is performed within a CAD software, the visualization of the resulting surface is modified but its underlying mathematical description remains unchanged. For a detailed review of trimming and related open challenges in IGA we refer to \citep{Marussig2018} and references therein.
In the subsequent Section, we introduce the basic mathematical setting in the scope of isogeometric methods defined on trimmed domains, following closely the notation used in \citep{Antolin2019,Buffa2019}.     

\subsection{Mathematical framework of trimming}
\begin{figure}[!ht]
	\centering
	\includegraphics[width=0.6\textwidth]{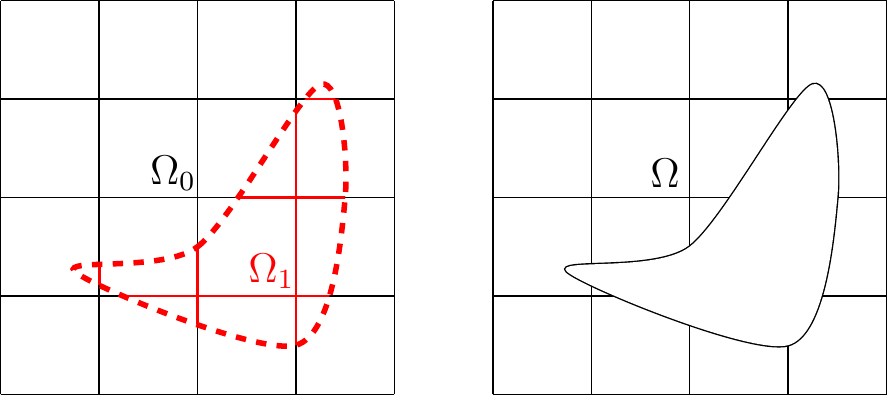}
	\caption{Example of a trimmed domain. From the original rectangular domain $\Omega_0$ the red domain $\Omega_1$ is trimmed away to produce the final domain $\Omega$. For a correct interpretation of the colors, the reader is referred to the online version of this manuscript.} 
	\label{fig:trimming_exaple}
\end{figure}

Let us define the domain $\Omega_0 \subset \mathbb{R}^d$, where $d$ is the dimension of the physical space of the problem at hand, described by a spline map $\mathbf{F}: \widehat{\Omega}_0 = [0, 1]^{d_r} \rightarrow \Omega_0$, where $\widehat{\Omega}_0$ denotes the parametric domain and $d_r$ its dimension. Let us also introduce $\Omega_1, \ldots, \Omega_N \subset \mathbb{R}^d$ Lipschitz-regular domains that define regions to be trimmed away from $\Omega_0$. Then, the physical domain can be written as follows:
\begin{align} \label{eq:domain}
\Omega = \Omega_0 \setminus \bigcup_{i = 1}^N \overline{\Omega}_i \, ,
\end{align}  
where an example is depicted in~\Cref{fig:trimming_exaple}. 
It is worth highlighting that the trimming operation does not change the underlying mathematical description of the original domain. Hence, elements and basis functions are constructed with respect to the non-trimmed domain $\Omega_0$. Given the definition of our computational domain in~\eqref{eq:domain}, we can now split its boundary in a part which coincides with the boundary of the original domain $\partial \Omega \cap \partial \Omega_0$ and a trimmed part $\partial \Omega \setminus \partial \Omega_0$. Furthermore, we denote with $\Gamma_D$ the Dirichlet part of boundary. In the following, if not specified otherwise, we will assume that $\Gamma_D \subset \partial \Omega \cap \partial \Omega_0$ and $\Gamma_D \neq \emptyset$, i.e. Dirichlet-type boundary conditions are not applied on the trimmed part of the boundary and, for simplicity, they are a non-empty set.

\subsection{Hierarchical B-spline basis} \label{subsec:hb}

Let us review the fundamental concepts behind the hierarchical B-spline basis, denoted by $\mathcal{HB}$ in the following. This allows us to build a basis that is locally refinable and therefore to overcome the limitations intrinsic to the tensor-product nature of B-splines and NURBS. Here, we summarize the definition of $\mathcal{HB}$ given in~\citep{Vuong2011,Kraft1997}.

\noindent Let $ V^0 \subset V^1 \subset \ldots \subset V^{N} $ be a sequence of $N+1$ nested multivariate splines spaces, determined by their corresponding degree and knot vectors. Then, the B-spline basis associated to the space $V^l$ of level $l$ is denoted by $\mathcal{B}^l := \left\lbrace b^{i,l} \, \vert \, i = 1, \ldots, N^l \right\rbrace$, where $N^l$ represents the dimension of $V^l$. Additionally, let us denote by $\mathcal{Q}^l_0$ the mesh associated to $\mathcal{B}^l$, where $Q \in \mathcal{Q}^l_0$ represents a cell of level~$l$. 

\noindent Let us also define the set $\boldsymbol{\Omega}^N_0 := \left\lbrace \Omega^0_0, \Omega^1_0, \ldots, \Omega^N_0 \right\rbrace$ as a hierarchy of subdomains of depth $N$ if the following holds:
\begin{align}\label{eq:subdomains}
\widehat{\Omega}_{0} = \Omega^0_0 \supset \Omega^1_0 \supset \ldots \supset \Omega^{N-1}_0 \supset \Omega^{N}_0 = \emptyset \, , 
\end{align}
and each subdomain $\Omega^l_0$ is the union of closed cells of level $l-1$, where $\Omega^0_0$ coincides with the initial parametric domain $\widehat{\Omega}_{0}$. The subscript 0 is used to remark that up to now we refer \textit{only} to the non-trimmed domain $\Omega_{0}$.
We can now define the hierarchical B-spline basis $\mathcal{HB}$. Given a sequence of spaces $\left\lbrace V^l \right\rbrace_{l=0,\ldots,N}$ (as defined above) with the corresponding B-spline bases $\left\lbrace \mathcal{B}^l \right\rbrace_{l=0,\ldots,N}$ and a hierarchy of subdomains $\boldsymbol{\Omega}^N_0$ of depth $N$, we define $\mathcal{HB}$ as the $(N-1)$-th step of the following recursive definition:
\begin{align*}
	\mathcal{HB}^0 = \, &\mathcal{B}^0 \\
	\mathcal{HB}^{l+1} = \, &\left\lbrace b \in \mathcal{HB}^l \, \vert \, \text{supp } b \not\subset \Omega^{l+1}_0 \right\rbrace \, \cup \\
	&\left\lbrace b \in \mathcal{B}^{l+1} \, \vert \, \text{supp } b \subset \Omega^{l+1}_0 \right\rbrace \, , \quad l=0,\ldots,N-2 \, .
\end{align*}
Consequently, we can introduce the (parametric) mesh of level $\widehat{\mathcal{Q}}_A^l$ and the resulting hierarchical mesh $\widehat{\mathcal{Q}}_0$ associated to $\mathcal{HB}$ as:
\begin{align*}
	\widehat{\mathcal{Q}}_{0,A}^l &= \left\lbrace Q \in \mathcal{Q}^l_0 \, \vert \, Q \subset \Omega^l_0 \wedge Q \not\subset \Omega^{l+1}_0 \right\rbrace \\
	\widehat{\mathcal{Q}}_0 &= \bigcup\limits_{l=0}^{N-1} \widehat{\mathcal{Q}}_{0,A}^l \, .
\end{align*}
Finally, in the following, we will refer to $b$ as an active function if $b \in \mathcal{HB}$ and as an active function of level $l$ if $b \in \mathcal{B}^l_A := \mathcal{HB} \cap \mathcal{B}^l$. An analogous terminology can be introduced for elements, where we denote $Q$ as an active cell if $Q \in \widehat{\mathcal{Q}}_0$ and as an active cell of level $l$ if $Q \in \widehat{\mathcal{Q}}_{0,A}^l$.
For further details we refer to~\citep{Garau2018}.

\subsection{Truncated Hierarchical B-spline basis} \label{subsec:thb}

We denote the truncated hierarchical B-spline basis as $\mathcal{THB}$, as introduced in~\citep{Giannelli2012}. The only difference between $\mathcal{HB}$ and $\mathcal{THB}$ is that in the latter the basis functions whose support overlaps finer elements are truncated, as described in~\citep{Giannelli2012}. This generates a basis that spans the same space as $ \mathcal{HB} $, but has in general better properties from a numerical standpoint, we refer to~\citep{Giannelli2013} for further details.
Let us first introduce the truncation operator. Namely, exploiting the fact that B-splines of level $l$ can be written as a linear combinations of B-splines of level $l+1$ with non-negative coefficients, the truncation operator with respect to level $l+1$ reads:  
\begin{align*}
\text{trunc}^{l+1}(b^{i,l}) = \sum \limits_{k=1}^{N^{l+1}} \tilde{c}^{\, k,l+1}(b^{i,l})  \, b^{k,l+1} \, ,
\end{align*}
where the coefficients are defined by:
\begin{align*}
\tilde{c}^{\, k,l+1}(b^{i,l}) = \begin{cases}
               0 \quad &\text{if } b^{k,l+1} \in \mathcal{HB}^{l+1} \cap \mathcal{B}^{l+1} \, ,\\
               c^{k,l+1} \quad &\text{otherwise,}
            \end{cases}
\end{align*}
where $c^{k,l+1}$ are the standard coefficients of the two-scale relation. Again, for further details we refer to~\citep{Garau2018} and references therein. Now, following~\citep{Giannelli2012}, we define the truncated hierarchical basis $\mathcal{THB}$ as the $(N-1)$-th step of the following recursive algorithm:
\begin{align*}
	\mathcal{THB}^0 = \, &\mathcal{B}^0 \\
	\mathcal{THB}^{l+1} = \, &\left\lbrace \text{trunc}^{l+1}(b) \, \vert \, b \in \mathcal{THB}^{l} \wedge \text{supp } b \not\subset \Omega^{l+1}_0 \right\rbrace \, \cup \\
	&\left\lbrace b \in \mathcal{B}^{l+1} \, \vert \, \text{supp } b \subset \Omega^{l+1}_0 \right\rbrace \, , \quad l=0,\ldots,N-2 \, .
\end{align*}
Finally, we can introduce the following discrete space:
\begin{align*}
X_h = \text{span }\left\lbrace b \circ \mathbf{F}^{-1} \, \vert \, b \in \mathcal{HB} \right\rbrace = \text{span }\left\lbrace b \circ \mathbf{F}^{-1} \, \vert \, b \in \mathcal{THB} \right\rbrace \, .
\end{align*}

\subsection{Truncated Hierarchical B-splines on trimmed domains}

Let us now trim the aforementioned domain $\Omega_0$ as defined in~\eqref{eq:domain} and let us make the assumption that $\Omega \subset \mathbb{R}^d$ is regular enough. 
In order to define a (truncated) hierarchical basis suitable for the analysis on $\Omega$, we follow the construction provided in~\citep[Section 4.5]{Hoellig_book}, which guarantees the linear independence of the basis.

Recalling the definition provided in~\eqref{eq:subdomains} for the non-trimmed case, we introduce a hierarchy of trimmed subdomains $\boldsymbol{\Omega}^N := \left\lbrace \Omega^0, \Omega^1, \ldots, \Omega^N \right\rbrace$ such that:
\begin{align}\label{eq:subdomains_trimmed}
\widehat{\Omega} = \Omega^0 \supset \Omega^1 \supset \ldots \supset \Omega^{N-1} \supset \Omega^{N} = \emptyset \, ,
\end{align}
where $\widehat{\Omega}$ denotes the trimmed parametric domain and $N$ represents the depth of the hierarchy. Notice that we dropped the subscript in our notation consistently with our definition of trimmed domain~\eqref{eq:domain}.
Then, let us define the support of a spline function $b$ restricted to the trimmed domain as:
\begin{align}\label{eq:supp_trimmed}
\text{supp}_{\, \widehat{\Omega}}(b) := \overline{\widehat{\Omega}} \cap \, \supp{b} \, .
\end{align}
Now, for every level $l$, we introduce the corresponding B-spline basis restricted to $\widehat{\Omega}$ as:
\begin{align}\label{eq:basis_trimmed}
\mathcal{B}^l_{\widehat{\Omega}} := \left\lbrace b \vert_{\widehat{\Omega}} \, : \, b \in \mathcal{B}^l \wedge \, \mathrm{meas}(\mathrm{supp}_{\, \widehat{\Omega}}(b) ) \neq 0 \right\rbrace \, .
\end{align}
With these definitions, the recursive algorithms for the construction of the $\mathcal{HB}$ and $\mathcal{THB}$ bases on trimmed domains are analogous to those introduced in~\Cref{subsec:thb,subsec:hb} in the non-trimmed case, where we need to replace the hierarchy of subdomains with its trimmed counterpart as defined in~\eqref{eq:subdomains_trimmed}, the support of functions with the definition in~\eqref{eq:supp_trimmed} and the B-spline basis of level $l$, $\mathcal{B}^l$, with its restriction to $\widehat{\Omega}$ as provided in~\eqref{eq:basis_trimmed}.
At this point, we highlight that it was proved in~\citep[Section 4.5]{Hoellig_book} that this construction guarantees the linear independence of the basis in the context of hierarchical B-splines.
We remark that the truncation operator does not affect this property (as it only acts reducing the original support) and therefore the same rationale can be applied to the $\mathcal{THB}$ basis, yielding again a basis that is suitable for the analysis in the scope of trimmed domains.
%
%

\noindent Then, we also modify the definition of parametric mesh of level and hierarchical mesh as follows:
\begin{align*}
	\widehat{\mathcal{Q}}_{A}^l &= \left\lbrace Q \in \widehat{\mathcal{Q}}_{0,A}^l \, : \, \mathrm{meas}( Q \cap \widehat{\Omega}) \neq 0 \right\rbrace \\
	\widehat{\mathcal{Q}} &= \left\lbrace Q \in \widehat{\mathcal{Q}}_0 \, : \, \mathrm{meas}( Q \cap \widehat{\Omega}) \neq 0 \right\rbrace \, ,
\end{align*}
from which the physical mesh reads:
\begin{align*}
	\mathcal{Q} = \left\lbrace \mathbf{F}(Q) \, : \, Q \in \widehat{\mathcal{Q}} \right\rbrace \, .
\end{align*}

\begin{remark}
We remark that the construction of hierarchical B-splines can be performed directly on the physical domain, replacing $b$ by $b \circ \mathbf{F}^{-1}$, and the subdomains $\Omega^l_0$ by $\mathbf{F}(\Omega^l_0)$. The same holds for the trimmed case.
\end{remark}

\begin{remark}
In order to provide an implementation of the (truncated) hierarchical basis as defined above, a proper function that checks the intersection between the support of basis functions and the trimmed domain $\Omega$ is crucial from an algorithmic standpoint.
However, this might not be straightforward to achieve into existing isogeometric codes. Here we provide an easier algorithm that identifies the active basis functions and retains the linear independence of the trimmed basis based on the standard implementation of (truncated) hierarchical B-splines. In particular, when a trimmed element $Q$ is marked for refinement, we also refine all those elements outside $\Omega$ (but inside $\Omega_0$) contained in the support of basis functions whose support incorporate $Q$. We refer to these elements as `ghost' cells. This procedure guarantees that when all the elements in the support of a function that intersects $\Omega$ are refined, those ones outside $\Omega$ will also be refined, and therefore the function will be deactivated. We summarize a possible implementation in Algorithm~\ref{alg:thb}.
\end{remark}

\begin{algorithm}[H] 
\begin{algorithmic}[1]
	\Procedure{Avoid$\_$linear$\_$dependence}{hierarchical mesh $\mathcal{Q}$} 
	\For{\textbf{each} level $l$ of $\mathcal{Q}$}
		\For{\textbf{each} $Q \in \mathcal{Q}^l$ marked for refinement}
				\State Get all functions with support on $Q$
				\State Get the supports of these functions				
				\State Get all ghost cells within these supports
				\State Mark the ghost cells for refinement 				
		\EndFor
	\EndFor	
	\State Refine marked elements
	\EndProcedure
\end{algorithmic} 
\caption{Algorithm for avoiding linear dependence of the $\mathcal{HB}$ or $\mathcal{THB}$ basis in the standard implementation.}\label{alg:thb}
\end{algorithm}

\noindent Lastly, we introduce the approximation space defined on the trimmed domain $\Omega$ as follows:
\begin{align*}
\widetilde{X}_h = \text{span }\left\lbrace  b \circ \mathbf{F}^{-1} \: \vert \, b \in \mathcal{HB} \right\rbrace = \text{span }\left\lbrace b \circ \mathbf{F}^{-1} \: \vert \, b \in \mathcal{THB} \right\rbrace \, ,
\end{align*}
where the construction of $\mathcal{HB}$ and $\mathcal{THB}$ takes trimming into account as explained in details above.
To avoid unnecessary complications at this stage, let us introduce a generic discrete variational problem of the form: find $u_h \in \widetilde{V}_h$ such that
\begin{align*}
	&a(u_h,v_h) = F(v_h)  \quad \forall v_h \in \widetilde{V}_h \, ,
\end{align*}
where $a(\cdot,\cdot)$ and $F(\cdot)$ are bilinear and linear forms, respectively, and they will be specified in the following according to the analyzed problem. Moreover, the choice of discrete space $\widetilde{V}_h \subset \widetilde{X}_h$ depends in general on the boundary conditions of such problem. We recall that in our setting boundary conditions are applied on $\partial \Omega \cap \partial \Omega_0$.

\section{A posteriori error estimator on trimmed domains} \label{sec:bubbles}

In the following Section, we introduce a variant of the family of error estimators studied in \citep{Bank1993} and we extend its isogeometric version, proposed in \citep{Juettler2010} and studied by the authors in~\citep{Antolin2019b} in the context of non-trimmed domains, to elliptic second and fourth order PDEs defined on trimmed domains. Then, we present a possible implementation of the proposed indicator which makes use of truncated hierarchical B-splines.

\subsection{The bubble error estimator}
The family of multi-level \textit{a posteriori} estimators was introduced in \citep{Bank1993} in the context of the $p$-version of the Finite Element Method and was first successfully applied to T-splines in \citep{Juettler2010} and later to hierarchical B-splines in~\citep{Vuong2011} in the scope of second-order elliptic PDEs. It has then been extended by the authors to Kirchhoff plates and Kirchhoff-Love shells in~\citep{Antolin2019b}. We now review the method and recast it into the framework of trimmed hierarchical IGA. Let us define the finite dimensional solution space $\widetilde{V}^p_h$ as the span of THB-splines basis functions of order $p$, as defined in \Cref{sec:trimming}.
Then, let us denote by $u_h \in \widetilde{V}^p_h$ the discrete solution of the problem at hand. 

\begin{figure}[!h]
	\begin{subfigure}[t]{0.475\textwidth}
	\centering
	\includegraphics[scale=1.0]{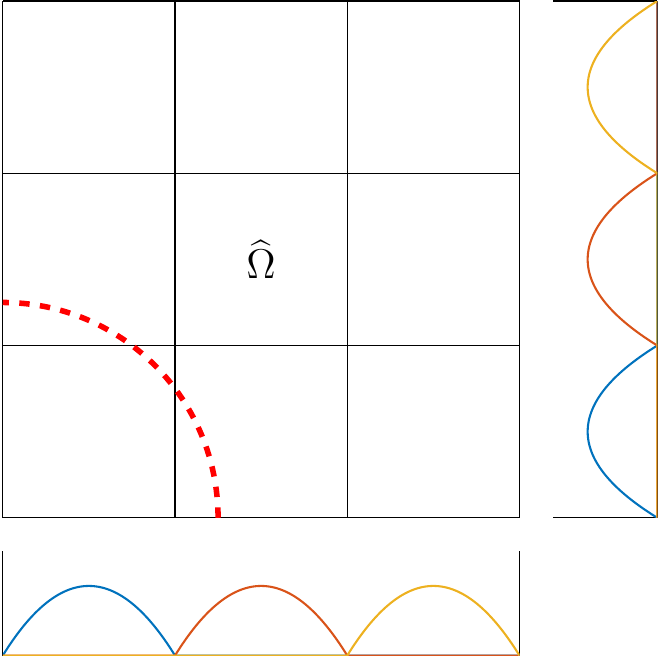}
	\caption{Bubble basis $p+1=2$ associated to an univariate knot vector $\Xi = [0 \, 0 \,0 \,1 \,1 \,2 \,2 \,3 \,3 \,3]$ in each parametric direction.}
	\label{fig:bubble_basis_p2}
	\end{subfigure}
	\begin{subfigure}[t]{0.475\textwidth}
	\centering
	\includegraphics[scale=1.0]{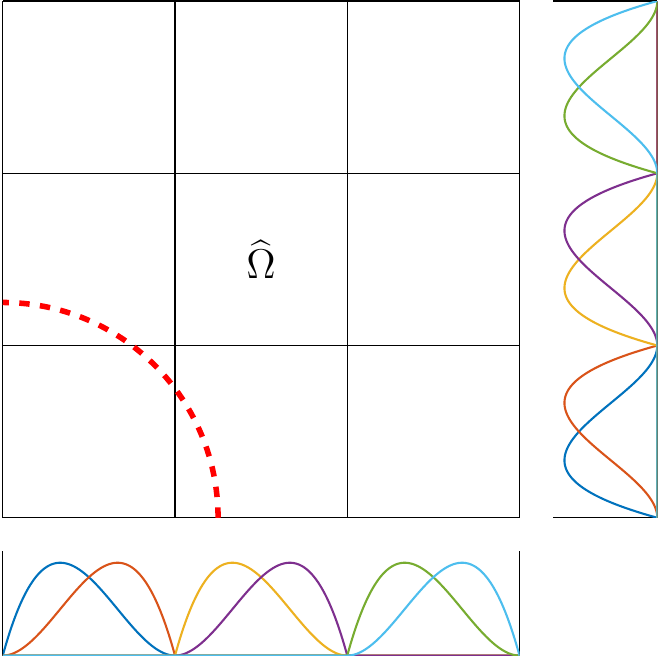}	
	\caption{Bubble basis $p+1=3$ associated to an univariate knot vector $\Xi = [0 \, 0 \, 0 \, 0 \, 1 \, 1 \, 1 \, 2 \, 2 \, 2 \, 3 \, 3 \, 3 \, 3]$ in each parametric direction.}
	\label{fig:bubble_basis_p3}
	\end{subfigure}		
	\caption{Example of construction of $C^0$ bubble functions on a trimmed parametric domain as tensor product of univariate functions for the cases $p+1=2,3$. The corresponding space is suitable for the error estimation of the solution of second order PDEs. We remark that we handle trimming at integration level by reparametrizing those elements that are cut.}
	\label{fig:bubble_basis}
\end{figure}

\noindent Now, starting from the degree and continuity of the aforementioned space $\widetilde{V}^p_h$ and for every level $l$, let us introduce the B-splines of degree $p+1$ and of reduced continuity, denoted as ${B}^{p+1,l}_{red}$, defined on $\mathcal{Q}_0^l$ (the rectangular grid of level $l$ introduced in~\Cref{sec:trimming}), where the chosen continuity depends on the problem at hand (see~\Cref{rem:bubble_continuity}). From an algorithmic point of view, this space is obtained by performing first degree elevation followed by increasing the multiplicity of the original knots. 
Now, from all the functions in ${B}^{p+1,l}_{red}$, we select only a suitable subset $\widetilde{B}^{p+1,l}_{act}$, which we characterize as:
\begin{align*}
\widetilde{B}^{p+1,l}_{act} = \lbrace b \in {B}^{p+1,l}_{red} \, \vert \text{ supp}(b) \subset Q \text{ for some } Q \in \widehat{\mathcal{Q}}_{A}^l \rbrace \, .
\end{align*}
In the following, we refer to these splines as \textit{bubble functions}, where an example defined on a single level is depicted in~\Cref{fig:bubble_basis}. We remark that by definition all functions in $\widetilde{B}^{p+1,l}_{act}$ have support on a single element. In fact, the bubble functions on each element always coincide with a subset of (scaled) Bernstein polynomials, and as a consequence they can be defined on a single reference element, which is then mapped to the elements of the mesh, analogously to basis functions in finite elements. We remark that this reference element is independent of the level. 

\noindent Then, let assume there exists a larger space $\widetilde{V}^p_h \subset \widetilde{W}^p_h \subset V$ such that the following decomposition holds:
\begin{align*}
\widetilde{W}^p_h = \widetilde{V}^p_h \oplus \widetilde{Z}^{p+1}_h \, ,
\end{align*} 
where $\widetilde{Z}^{p+1}_h$ is the space (defined over a trimmed geometry) in which we seek a good estimation of the error $e_h \approx e = \left\lVert u - u_h \right\rVert$, in a suitable norm $\left\lVert \cdot \right\rVert$. In particular, we can characterize $\widetilde{Z}^{p+1}_h$ in a multi-level fashion as follows:
\begin{align*}
\widetilde{Z}^{p+1}_{h} = \bigcup_{l=0}^N \widetilde{Z}^{p+1}_{h,l} \, ,
\end{align*}
where 
\begin{align*}
\widetilde{Z}^{p+1}_{h,l} = \text{span }\lbrace b \circ \mathbf{F}^{-1} \vert \, b \in \widetilde{B}^{p+1,l}_{act} \rbrace \, .
\end{align*}
Namely, $\widetilde{Z}^{p+1}_{h,l}$ is the space spanned by active B-splines of level $l$ obtained by degree elevation and knot insertion as discussed above, such that their support is compact and overlaps exactly with one active element (trimmed or non-trimmed) $Q$ of level $l$, where we postpone the discussion on the required continuity of $b$ to~\Cref{rem:bubble_continuity}. We highlight that, in the trimmed case, we must replace the standard definition of support of a function with its trimmed counterpart, which was previously denoted as $\mathrm{supp}_{\, \widehat{\Omega}}(\cdot)$.
We are now ready to define the a posteriori error estimate $e_h$ as the solution to the following problem: find $e_h  \in \: \widetilde{Z}^{p+1}_h$ such that
\begin{subequations}\label{eq:bubbleProblem}
 \begin{equation}\label{eq:bubbleSystem}
 a(e_h,b_h) = F(b_h) - a(u_h,b_h)  \qquad \forall b_h \in \widetilde{Z}^{p+1}_h \, .
 \end{equation}  
Notice that due to the compact support property of $b_h$, the linear system corresponding to the discrete error weak form \eqref{eq:bubbleSystem} is block diagonal, where each block corresponds to a single element, and its size is given by the number of bubble functions. Lastly, we compute the element-wise error estimator $\eta_{Q}$ as:
\begin{equation}\label{eq:bubbleSystem_bis}
\eta_{Q} = C_a \left\lVert e_h \right\rVert_{E(Q)}  \quad \forall Q \in \mathcal{Q} \, ,
\end{equation}  
\end{subequations}
where $\left\lVert \cdot \right\rVert_{E(Q)}$ denotes the energy norm restricted to an element $Q$ of the hierarchical mesh $\mathcal{Q}$. 
Additionally, we introduce the heuristic constant $C_a$ which, considering all our numerical experiments, seems to be independent from the chosen degree and from the problem at hand. In the scope of this work, we set $C_a = 3$.

It is worth highlighting that up to this point the derivation is completely independent from the problem we are solving, in the sense that once bilinear and linear forms and a suitable bubble space are chosen, the methodology follows the same aforementioned steps. Moreover, trimming can be treated naturally in this framework.
This confirms the simplicity of the proposed estimator on trimmed domains and that the method can be implemented into existing isogeometric codes in a straightforward manner, e.g. see the pseudo-code summarized in Algorithm~\ref{alg:estimate}.
Furthermore, our technique is computationally cheap and embarrassingly parallelizable due to the definition of the disjoint bubble space $\widetilde{Z}^{p+1}_h$ and the choice of solving \eqref{eq:bubbleProblem} level-wise over the hierarchical mesh $\mathcal{Q}$. 

\begin{algorithm}[H] 
\begin{algorithmic}[1]
	\Procedure{Estimate$\_$error}{numerical solution $u_h$, hierarchical mesh $\mathcal{Q}$} 
	\State Initialize vector $\eta$
	\For{\textbf{each} level $l$ of $\mathcal{Q}$}
		\State Build bubble space $\widetilde{Z}^{p+1}_{h,l}$
		\For{\textbf{each} $Q \in \mathcal{Q}^l$}
			\State Given $u_h$, solve locally on element $Q$ the additional system \eqref{eq:bubbleSystem} 
			\State Compute the element-wise indicator $\eta_{Q}$ \eqref{eq:bubbleSystem_bis}
			\State Store $\eta_{Q}$ in $\eta$
		\EndFor
	\EndFor	
	\State Return $\eta$
	\EndProcedure
\end{algorithmic} 
\caption{Bubble error estimator algorithm}\label{alg:estimate}
\end{algorithm}

\begin{remark}[On the choice of the bubble space]\label{rem:bubble_continuity}
We note that the choice of the bubble space $\widetilde{Z}^{p+1}_h$ has to be compatible with the underlying bilinear form. Therefore, for every level $l$, we select the continuity of ${B}^{p+1,l}_{red}$ such that all functions are at least $C^0$-continuous globally for second order problems (such as linear elasticity) and at least $C^1$-continuous for fourth order problems (such as the Kirchhoff-Love shell). Finally, we remark that when inhomogeneous Neumann boundary conditions are imposed, additional functions should be taken to correctly capture the error contribution on the boundary, for further details see~\citep{Antolin2019b}.
\end{remark}

\begin{remark}[Mark and refine]
Following the standard procedure in adaptivity, once the element-wise error $\eta_{Q}$ has been computed for all elements of the (trimmed) hierarchical mesh $\mathcal{Q}$, we mark elements for refinement by employing the so-called \textit{maximum strategy}, for additional details we refer to~\citep{Ainsworth1997}. In brief, given a user-defined threshold $\gamma \in (0,1)$, all elements such that
\begin{align*}
\eta_{Q} > \gamma \, \tilde{\eta}_{Q} \, , \quad \text{where} \quad \tilde{\eta}_{Q} = \max_{Q \in \mathcal{Q}} \, \eta_{Q} \, ,
\end{align*}
are selected for refinement. Moreover, before actually performing refinement of the corresponding $\mathcal{THB}$ basis, our algorithm is designed to preserve the admissibility (as defined in~\citep{Buffa2016a}) of the hierarchical mesh $\mathcal{Q}$ between consecutive iterations of the adaptive procedure. 
From a practical point of view, this means that we prevent the interaction between functions belonging to very fine and very coarse levels. In our numerical examples, we will set the class of admissibility $m$ to be $m = p$ and $m = p-1$ for second order and fourth order PDEs, respectively. We refer the reader to~\citep{Buffa2016a,Buffa2017,Bracco2018} for a detailed description of the concept of admissibility and its application to adaptive refinement.
\end{remark}

\newcommand{\graphDir}{pictures/numericalExamples/graphs}
\newcommand{\dataDir}{pictures/numericalExamples/data}

\section{Numerical results} \label{sec:numerical_examples}

In the following Section we present several numerical experiments for various linear elliptic second- and fourth-order PDEs that assess the good performance of the proposed error estimator in steering adaptive simulations on trimmed geometries. All results have been obtained with the open-source and free Octave/Matlab isogeometric package \textit{GeoPDEs} \citep{Vazquez2016} in combination with the re-parametrization tool studied in~\citep{Antolin2019}.
We remark that, unless stated otherwise, we set the marking parameter $\gamma = 0.5$ in all our examples. Lastly, to reduce the detrimental effects of trimming on the conditioning number, we apply a simple diagonal scaling, which can be seen as a Jacobi preconditioner. This seems to suffice in all the following numerical experiments. For a more detailed discussion of the conditioning issues stemming from trimming we refer to~\citep{Deprenter2017,Deprenter2019}.


\subsection{The Poisson problem}

We start our numerical investigation from the well-known Poisson equation. Let us briefly recall the strong form of the problem at hand:
\begin{alignat}{2}\label{eq:strongFormLaplace}
	- \varDelta u &= f \quad &&\text{in} \quad \Omega \\ 
	u &= \tilde{u} \quad &&\text{on} \quad \Gamma_D \nonumber \\ 
	\frac{\partial u}{\partial n} &= \tilde{g} \quad &&\text{on} \quad \Gamma_N \nonumber \, ,
\end{alignat} 
where $\Gamma_D \subset \partial \Omega \cap \partial \Omega_0$ and $\Gamma_N$ denote the Dirichlet and Neumann part of the boundary $\partial \Omega$, respectively, and it holds $\overline{\Gamma}_D \cup \overline{\Gamma}_N = \partial \Omega$ and $\Gamma_D \cap \Gamma_N = \emptyset$. Additionally, $f \in L^2(\Omega)$ represents the given source term, $\boldsymbol{n}$ the outward normal to the boundary, and $\tilde{u} \in H^{\frac{1}{2}}(\Gamma_D)$ and $\tilde{g} \in H^{-\frac{1}{2}}(\Gamma_N)$ the prescribed Dirichlet and Neumann data, respectively. We will detail the definition of these quantities according to the proposed example. 
Following the standard derivation, we can write the discrete weak formulation of \eqref{eq:strongFormLaplace} as follows: find $u_h  \in \: \widetilde{V}_h$ such that
\begin{align}\label{eq:weakFormLaplaceDiscrete}
 &a(u_h,v_h) = F(v_h)  \qquad \forall v_h \in \widetilde{V}_h \, ,
\end{align} 
where $\widetilde{V}_h$ is the finite-dimensional subspace defined by the (truncated hierarchical) B-spline basis.
The bilinear form $a(\cdot,\cdot)$ and the linear functional $F(\cdot)$ can be expanded as:
\begin{align*}
 a(u_h,v_h) &= \int_{\Omega} \nabla v_h \, \cdot \, \nabla u_h \, \text{d}\Omega   \\
 F(v_h) &= \int_{\Omega} f v_h \, \text{d}\Omega + \int_{\Gamma_N} \tilde{g} \, v_h \, \text{d}\Gamma \, , 
\end{align*}  
respectively.
\color{black}

\subsubsection{Singular problem, line singularity}  
\begin{figure}
	\centering
		\centering
			\includegraphics[scale=0.65]{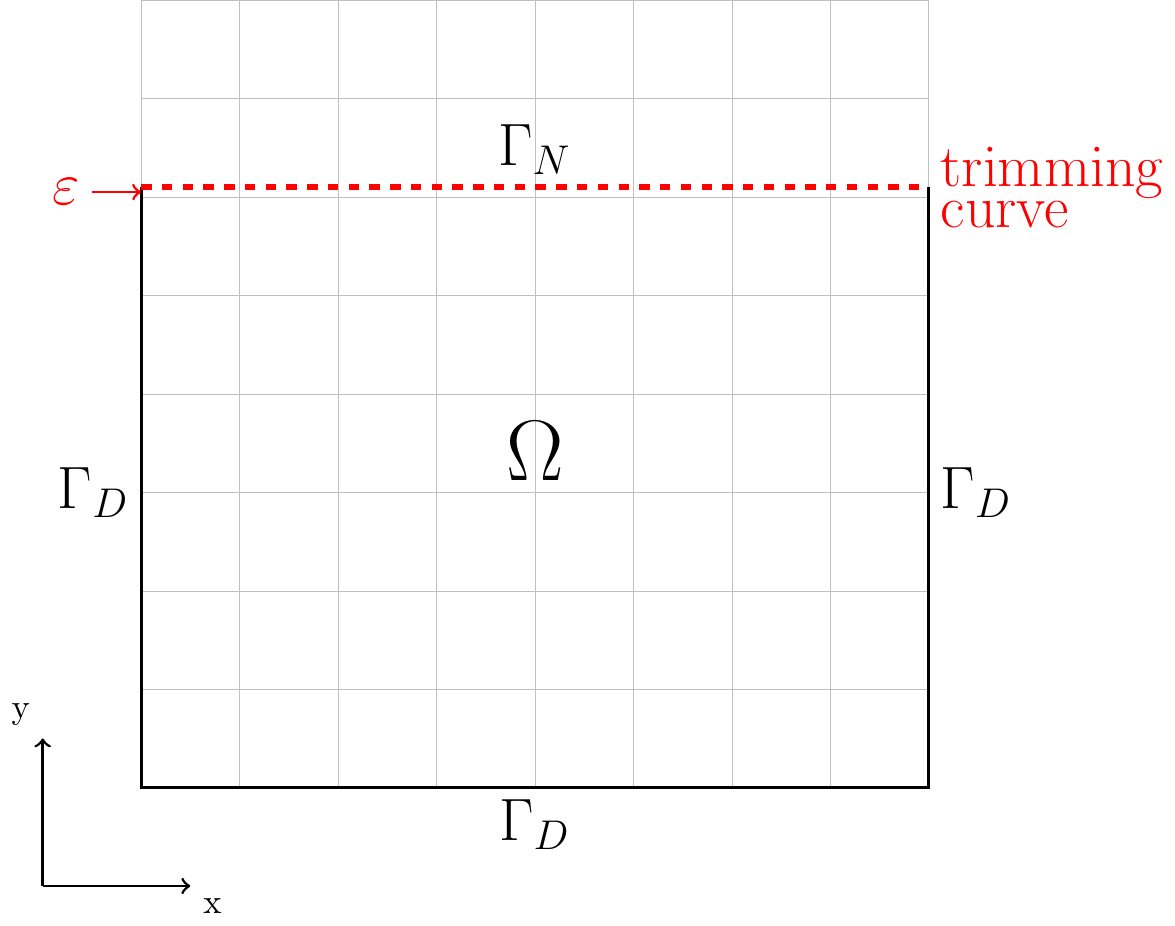}
		\caption{Problem description for the Poisson singular problem.} \label{fig:trimmed_square_geo}
\end{figure}

\begin{figure}
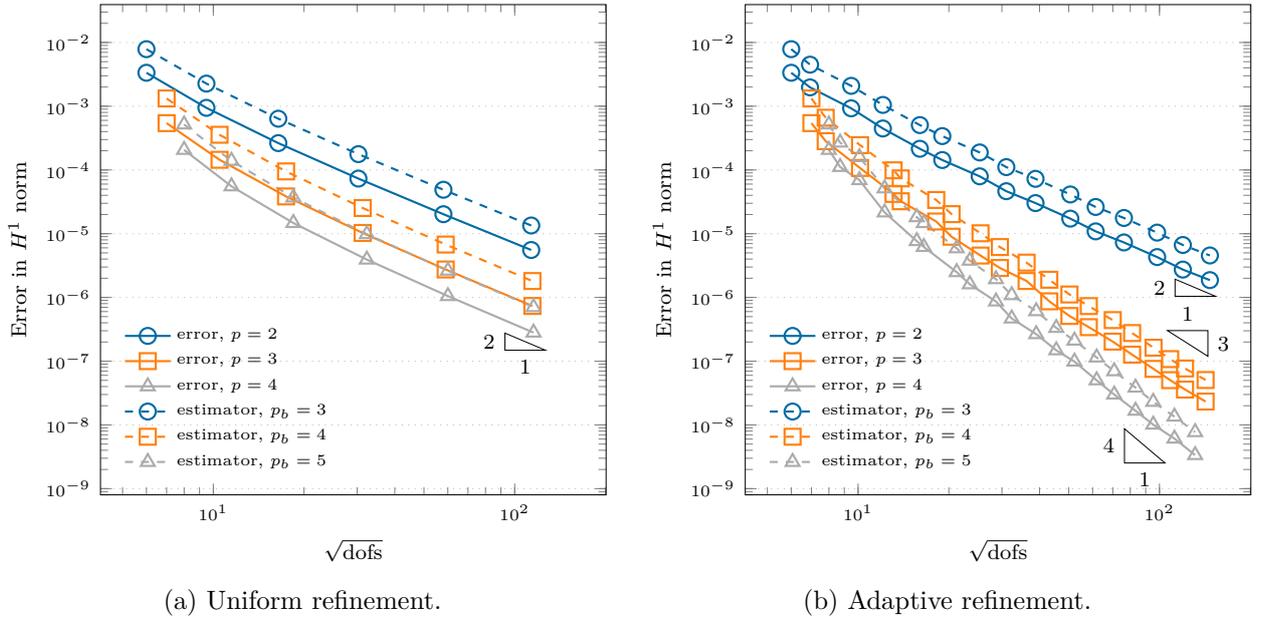

	\centering
	\begin{subfigure}[t]{0.495\textwidth}
		\centering
		\input{\graphDir/testcase_laplacian_2D_line_singularity_eps_1e5_trimmed_GR.tex}
		\caption{Uniform refinement.}\label{fig:convergence_lineSingularity_GR}
	\end{subfigure}
	\hfill
	\begin{subfigure}[t]{0.495\textwidth}
		\centering
		\input{\graphDir/testcase_laplacian_2D_line_singularity_eps_1e5_trimmed_MS.tex}
		\caption{Adaptive refinement.}\label{fig:convergence_lineSingularity_MS}
	\end{subfigure}
	\caption{Study of the convergence of the error measured in the $H^1$-norm and the bubble estimator for the Poisson singular problem. Comparison of uniform and adaptive refinements.}
\label{fig:convergence_lineSingularity}
\end{figure}

In this problem, we consider a trimmed computational domain defined by $\Omega = [0,1] \times [0,\frac{3}{4} + \varepsilon]$, as depicted in~\Cref{fig:trimmed_square_geo}, where $\varepsilon = 10^{-5}$ is chosen. As exact solution, we take a singular function of the form $u = x^\alpha y^\alpha$ with $\alpha = 2.4$. The imposed boundary conditions are also given in~\Cref{fig:trimmed_square_geo} and are computed such that they fulfill the exact solution, as well as the given source term $f$.
We remark that with this choice of parameter it holds $u \in H^2(\Omega) \setminus H^3(\Omega)$. From classical \textit{a priori} error analysis results \citep{Bazilevs2006}, it is known that this reduction in regularity hinders the rate of convergence achievable with uniform refinement of the mesh. In particular, the estimate reads as follows:
\begin{align*}
\left\lVert u - u_h \right\rVert_{H^1(\Omega)} \leq C h^{m - 1} \left\lVert u \right\rVert_{H^{r}(\Omega)} \qquad r > 1 \, ,
\end{align*}
where $r$ represents the regularity of the exact solution and $m$ is defined as $\min{(r,p+1)}$. For our choice of $\alpha$, optimal rates of convergence are achieved in case of uniform refinement only for bi-quadratic B-splines whereas higher degree discretizations suffer from the lack of regularity, see~\Cref{fig:convergence_lineSingularity_GR}.    
As shown in~\Cref{fig:convergence_lineSingularity_MS}, it is clear that optimal rates of convergence are recovered for all degrees $p=2,3,4$ with an adaptive simulation driven by the proposed estimator. Moreover, this example highlights the benefits of using an adaptive strategy in terms of accuracy per degree-of-freedom.  

\subsection{Linear elasticity}

In the following, we want to analyze the well-known problem of elastostatics, for more details we refer to \citep{Hughes2000b}. We consider a two-dimensional body $\Omega \subset \mathbb{R}^2$ where again we assume the boundary $\partial \Omega$ of such domain to consist of two disjoint parts $\Gamma_D $ and $\Gamma_N $, representing the Dirichlet and Neumann part of the boundary, respectively. We also assume that $\Gamma_D \subset \partial \Omega \cap \partial \Omega_0$.
For completeness, let us recall the strong formulation of the problem at hand: 
\begin{alignat}{2}\label{eq:strongFormElasticity}
	 - \nabla \cdot \boldsymbol{\sigma}(\boldsymbol{u}) &= \boldsymbol{f} \quad &&\text{in} \quad \Omega \\ 
	 \boldsymbol{\sigma}(\boldsymbol{u}) &= \mathbb{C} : \boldsymbol{\varepsilon}(\boldsymbol{u}) \quad &&\text{in} \quad \Omega \nonumber \\  
	 \boldsymbol{\varepsilon}(\boldsymbol{u}) &= \dfrac{1}{2} \left( \nabla \boldsymbol{u} + (\nabla \boldsymbol{u})^T \right) \quad &&\text{in} \quad \Omega \nonumber \\ 
	 \boldsymbol{\sigma}(\boldsymbol{u}) \cdot \boldsymbol{n} &= \boldsymbol{\tilde{g}} \quad &&\text{on} \quad \Gamma_N \nonumber \\  
	 \boldsymbol{u} &= \boldsymbol{\tilde{u}} \quad &&\text{on} \quad \Gamma_D \nonumber \, ,
\end{alignat}
where $\boldsymbol{u}$ represents the sought displacement of the structure and $\boldsymbol{\sigma}$, $\boldsymbol{f}$, $\mathbb{C}$, $\boldsymbol{\varepsilon}$, $\boldsymbol{\tilde{g}}$, $\boldsymbol{\tilde{u}}$, denote the stress tensor, body force, material tensor, infinitesimal strain tensor (defined as the symmetric part of the gradient of $\boldsymbol{u}$), prescribed traction and prescribed displacement, respectively. In the scope of this work, we will consider the body to be homogeneous and isotropic. 
Similarly to the previous case of the Poisson equation, 
the discrete weak formulation of \eqref{eq:strongFormElasticity} reads: find $\boldsymbol{u_h}  \in \: \widetilde{V}_h$ such that
\begin{align}\label{eq:weakFormElasticityDiscrete}
 &a(\boldsymbol{u_h},\boldsymbol{v_h}) = F(\boldsymbol{v_h})  \qquad \forall \boldsymbol{v_h} \in \widetilde{V}_h \, ,
\end{align} 
where the bilinear form $a(\cdot,\cdot)$ and the linear functional $F(\cdot)$ are given as:
\begin{align*}
 a(\boldsymbol{u_h},\boldsymbol{v_h}) &= \int_{\Omega} \boldsymbol{\sigma}(\boldsymbol{u_h}) \, : \, \boldsymbol{\varepsilon}(\boldsymbol{v_h}) \, \text{d}\Omega  \\ 
 F(\boldsymbol{v_h}) &= \int_{\Omega} \boldsymbol{f} \, \cdot \, \boldsymbol{v_h} \, \text{d}\Omega + \int_{\Gamma_N} \boldsymbol{\tilde{g}} \, \cdot \, \boldsymbol{v_h} \, \text{d}\Gamma \, .
\end{align*}  

\color{black}

\subsubsection{Infinite plate with a hole}
\begin{figure}[!ht]
	\centering
			\includegraphics[scale=0.9]{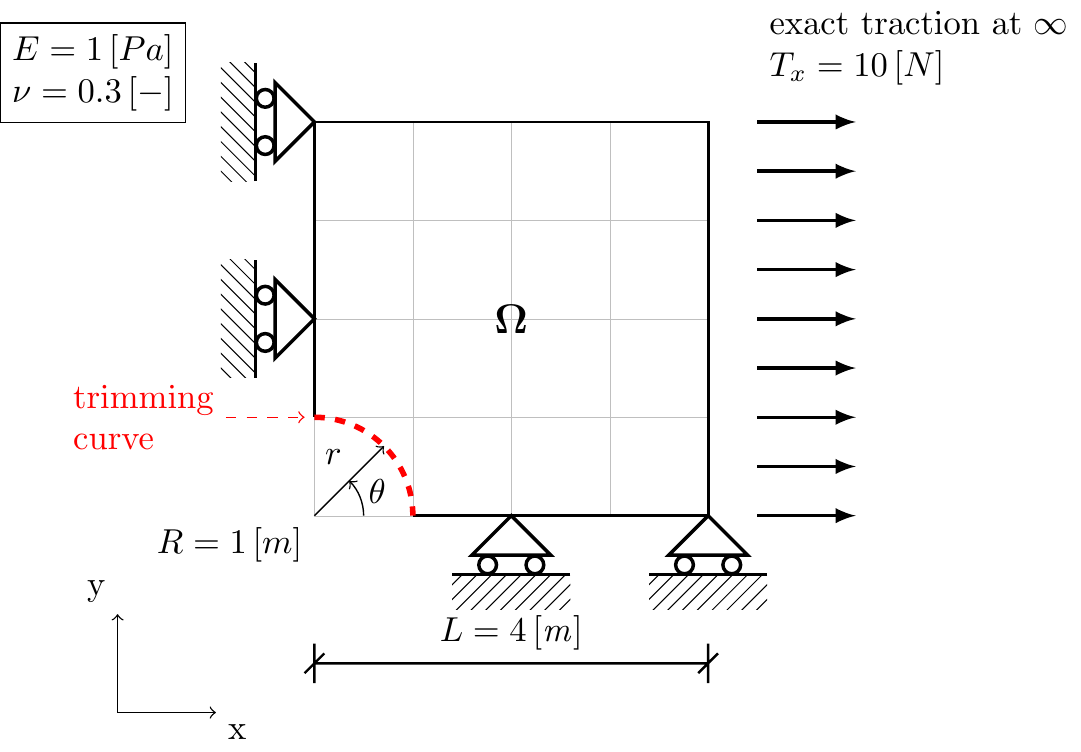}
	\caption{Problem setup and boundary conditions for the plate with a hole benchmark.}
 \label{fig:plate_with_hole_geo}
\end{figure}

\begin{figure}[!ht]
	\centering
	\begin{subfigure}[t]{0.495\textwidth}
	\input{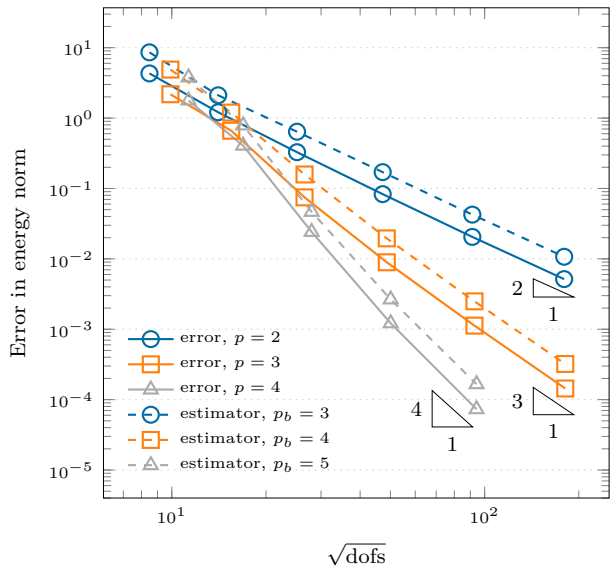}
		\caption{Uniform refinement.}
	\end{subfigure}
	\hfill
	\begin{subfigure}[t]{0.495\textwidth}
	\input{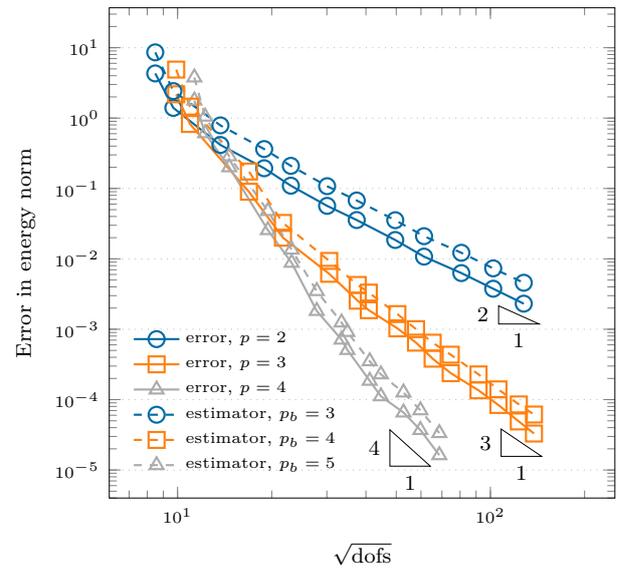}
		\caption{Bubbles estimator.}
	\end{subfigure}	
	\caption{Study of the convergence of the error measured in the energy norm and the bubble estimator for the plate with a hole. Comparison of uniform and adaptive refinements.}
\label{fig:convergence_plate_with_hole}
\end{figure}

\begin{figure}[!ht]
		\centering
		\begin{subfigure}[t]{0.495\textwidth}
			\centering
			\includegraphics[width=0.75\textwidth]{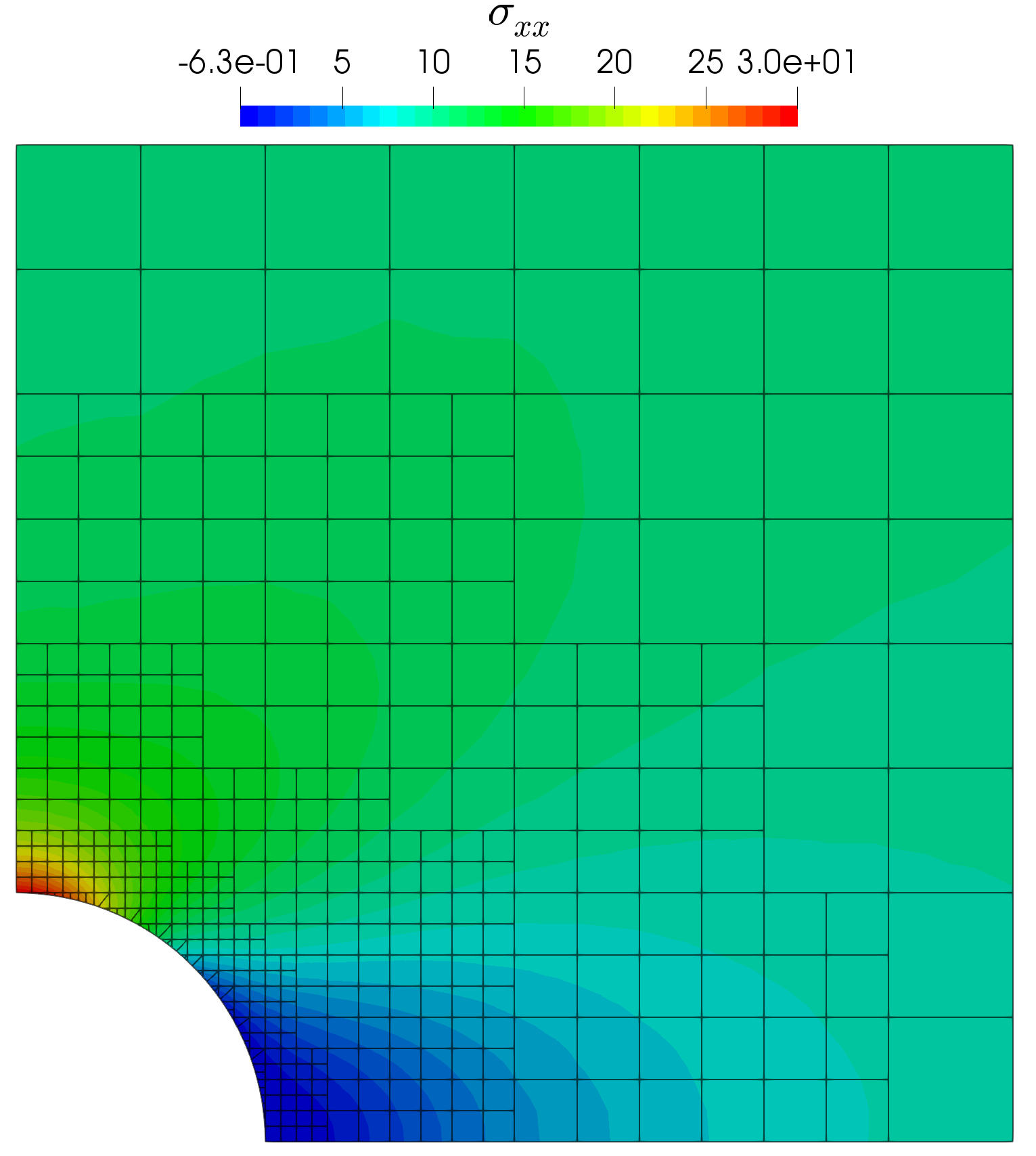}
			\caption{Mesh at iteration $k = 4$.}
		\end{subfigure}
		\begin{subfigure}[t]{0.495\textwidth}
			\centering
			\includegraphics[width=0.75\textwidth]{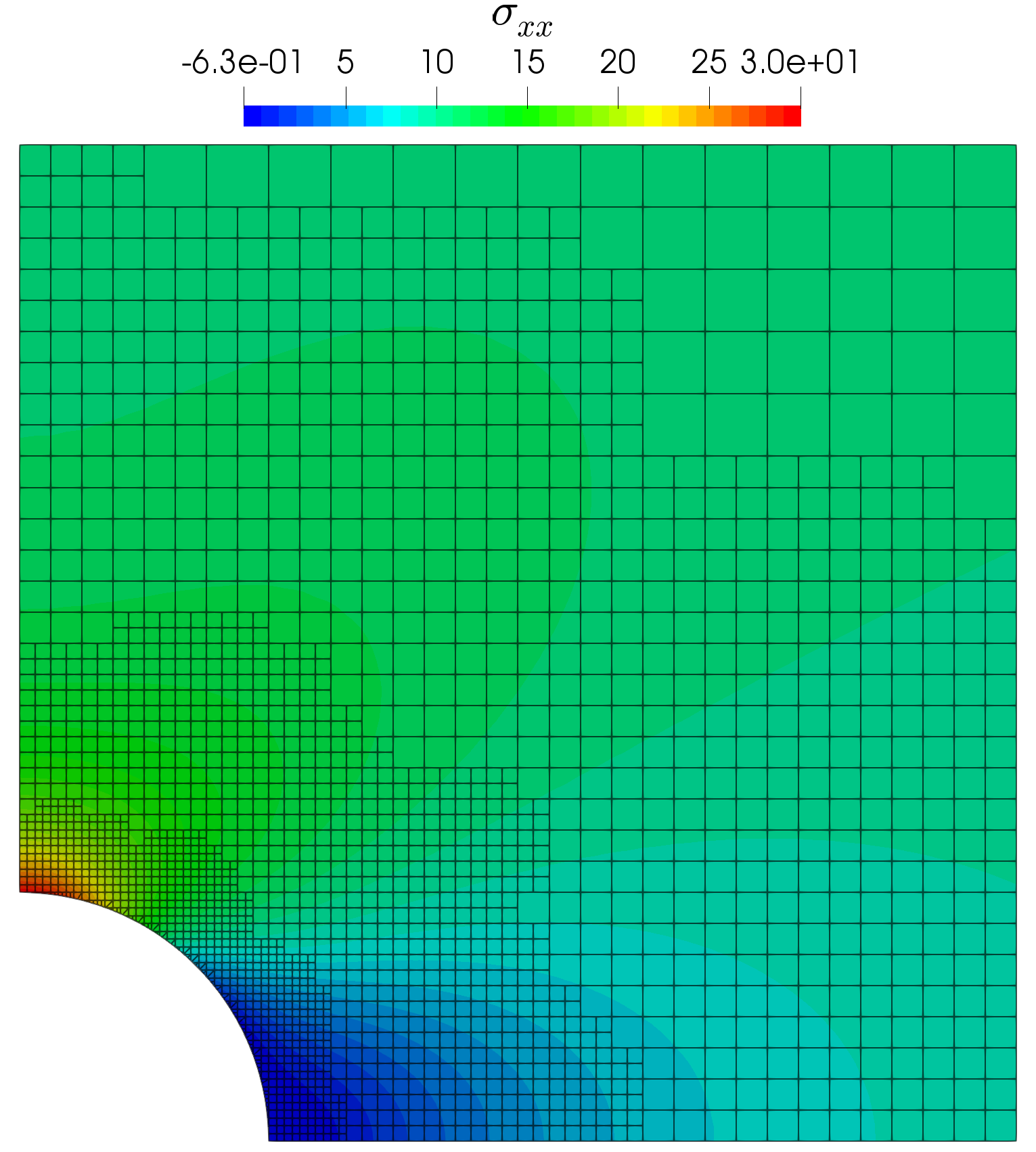}
			\caption{Mesh at iteration $k = 8$.}
		\end{subfigure}
		\begin{subfigure}[t]{0.495\textwidth}
			\centering
			\includegraphics[width=0.75\textwidth]{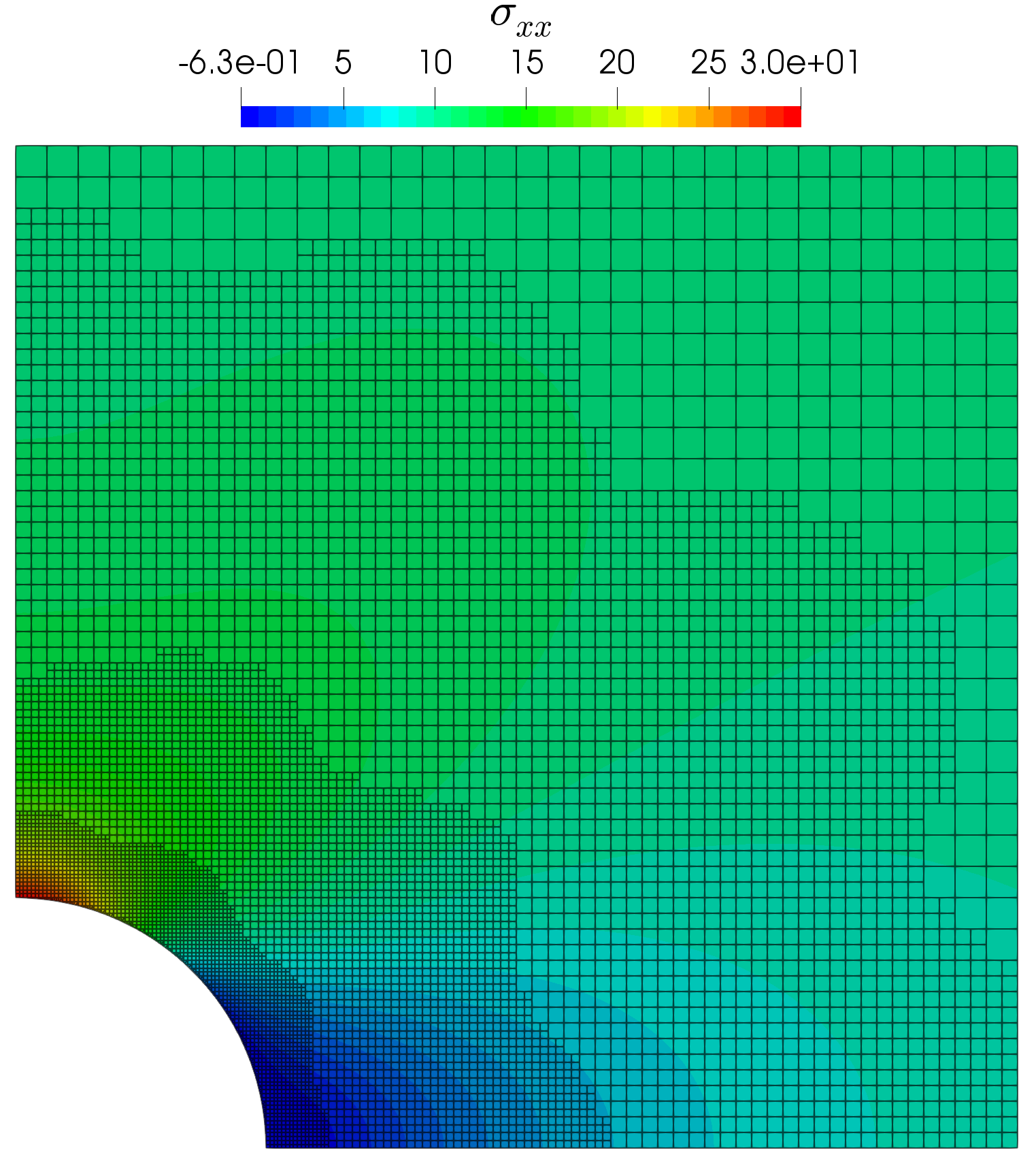}
			\caption{Mesh at iteration $k = 11$.}
		\end{subfigure}
		\begin{subfigure}[t]{0.495\textwidth}
			\centering
			\includegraphics[width=0.75\textwidth]{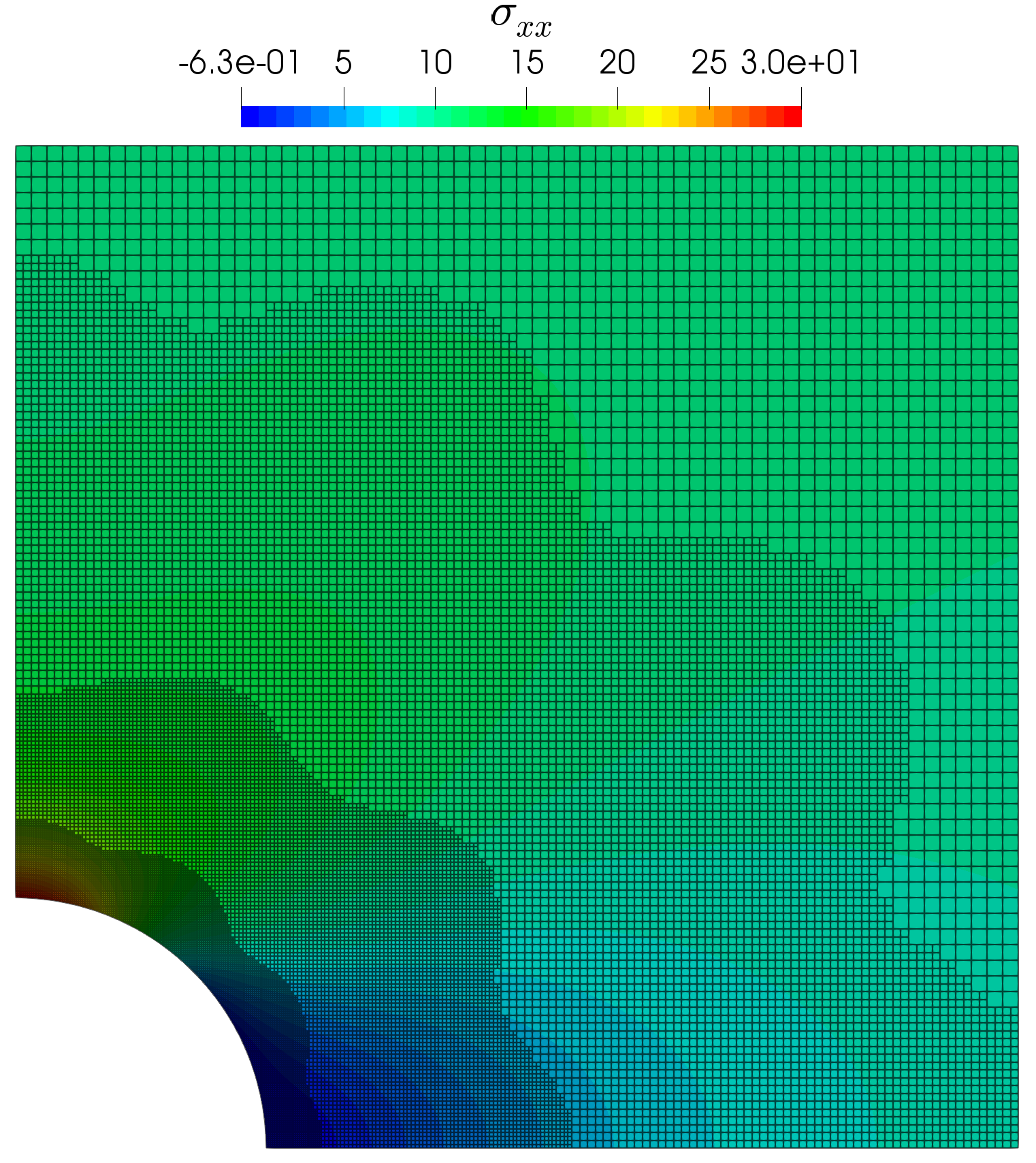}
			\caption{Mesh at iteration $k = 14$.}
		\end{subfigure}
			\caption{Mesh and stress $\sigma_{xx}$ at different steps of the adaptive loop driven by the bubble error estimator for the plate with a hole benchmark, solution obtained using THB-splines of degree $p=2$. Results obtained by setting the marking parameter to $\gamma = 0.55$.}
		\label{fig:mesh_plate_with_hole}
\end{figure}

In the next example, we analyze the well-known benchmark of an infinite plate with a circular hole. The geometry and problem setup are depicted in~\Cref{fig:plate_with_hole_geo}, where thanks to symmetry we model only one quarter of the plate. In particular, we use a Cartesian mesh trimmed by a curve corresponding to a quarter of a circle (light grey grid and dashed, red curve in~\Cref{fig:plate_with_hole_geo}, respectively).
We recall that the solution $\boldsymbol{u}$ is given as a function of the polar coordinates $(r,\theta)$ and reads \citep{Gould1999}:
\begin{align*}
\boldsymbol{u}(r,\theta) = \frac{T_x R}{8 \mu}
\begin{pmatrix}
\left(\frac{r}{R}(\kappa+1)\cos(\theta) + 2 \frac{R}{r}((1+\kappa)\cos(\theta)+\cos(3\theta))-2 \frac{R^3}{r^3}  \cos(3 \theta) \right)\\
\left(\frac{r}{R}(\kappa-3)\sin(\theta) + 2 \frac{R}{r}((1-\kappa)\sin(\theta)+\sin(3\theta))-2 \frac{R^3}{r^3}  \sin(3 \theta) \right)
\end{pmatrix} \, ,
\end{align*}
where $\kappa = 3-4\nu $ is the so-called Kolosov constant for the plane strain case, $\mu$ represents the second Lam\`{e} parameter, $T_x$ is the exact traction applied at infinity and $R$ denotes the radius of the hole. 
The solution $\boldsymbol{u}$ is smooth and therefore we expect optimal rates of convergence for both uniform and adaptive refinements. However, due to the presence of a hole, a stress concentration appears in the proximity of the bottom left corner. This feature is correctly detected and resolved by the error estimator, which yields a faster convergence in the pre-asymptotic regime of the adaptive simulation, as shown in~\Cref{fig:convergence_plate_with_hole}. Finally, in~\Cref{fig:mesh_plate_with_hole}, the hierarchical meshes obtained at various iterations $k=4,8,11,14$ of the adaptive loop are presented, where we remark that, in the visualization, the triangular elements close to the trimming curve are those provided by the re-parametrization tool for integration~\citep{Antolin2019}. Here, we notice how the refinement is at first steered around the hole where the stress concentration is located and then it gradually propagates into the rest of the computational domain. 

\FloatBarrier

\subsection{Kirchhoff-Love shell}

Next, we test the performance of the proposed method in the context of trimmed Kirchhoff-Love shells. 
The discrete weak formulation of the Kirchhoff-Love problem can be summarized as follows \citep{Kiendl2009,Cirak2006}: find $\boldsymbol{u_h} \in \widetilde{V}_h$ such that
\begin{align}\label{eq:modelWeakForm}
a\left(\boldsymbol{u_h} ,\boldsymbol{v_h} \right) = F \left(\boldsymbol{v_h} \right)  \qquad \forall \boldsymbol{v_h} \in \widetilde{V}_h \, ,
\end{align}
where $a(\cdot,\cdot)$ is a continuous and strongly coercive bilinear form and $F(\cdot)$ is a continuous linear functional. They can be written, respectively, as:
\begin{subequations}\label{eq:weak_KL}
\begin{align*}
	a(\boldsymbol{u_h},\boldsymbol{v_h}) &= \int_{\Omega} \boldsymbol{\varepsilon}(\boldsymbol{v_h}) \, : \, \boldsymbol{n}(\boldsymbol{u_h}) \, \text{d}\Omega + \int_{\Omega} \boldsymbol{\kappa}(\boldsymbol{v_h}) \, : \, \boldsymbol{m}(\boldsymbol{u_h}) \, \text{d}\Omega \, ,\\
	F(\boldsymbol{v_h}) &= \int_{\Omega} \boldsymbol{v_h} \cdot \boldsymbol{b} \, \text{d}\Omega + \int_{\Gamma_N} \boldsymbol{v_h} \cdot \boldsymbol{p} + \boldsymbol{\omega}(\boldsymbol{v_h}) \cdot \boldsymbol{r} \, \text{d}\Gamma  \, . 
\end{align*}
\end{subequations}
where $\boldsymbol{\varepsilon}$, $\boldsymbol{\kappa}$ denote the membrane and bending strain tensors, respectively, and $\boldsymbol{n}$, $\boldsymbol{m}$ are their energetically conjugate stress resultants. Additionally, $\boldsymbol{b}$ is the given body load, and $\boldsymbol{p}$ and $\boldsymbol{r}$ are the prescribed boundary forces and moments, respectively. Similarly to the examples above, we assume that $\Gamma_D \subset \partial \Omega \cap \partial \Omega_0$. For further details we refer to~\citep{Kiendl2009,Cirak2006}.

\color{black}

\subsubsection{Test case 1: Scordelis-Lo roof with an elliptic hole}
\begin{figure}[!ht]
	\centering
	\includegraphics[width=0.85\textwidth]{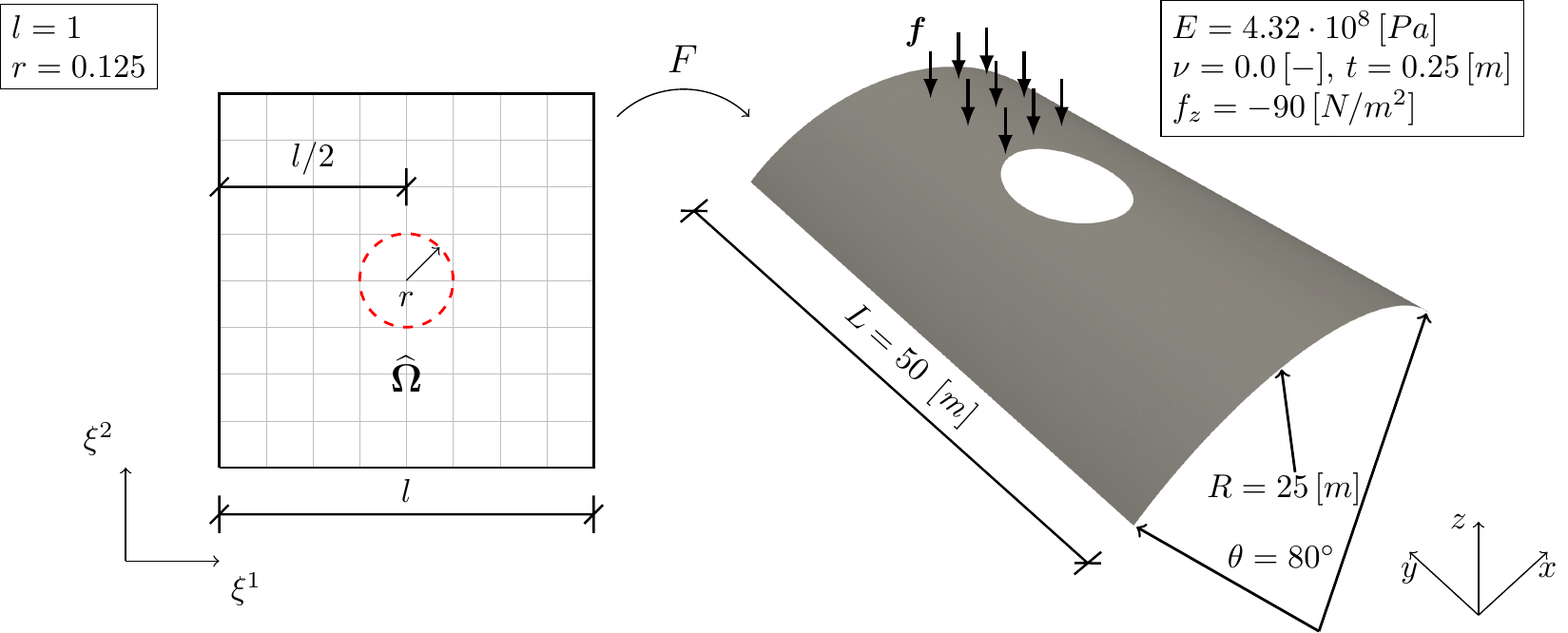}
	\caption{Geometry and physical parameters of the Scordelis-Lo roof with an elliptic hole (only part of the gravity load is shown for visualization purposes).} 
	\label{fig:scordelis_with_hole_geo}
\end{figure}

\begin{figure}[!ht]
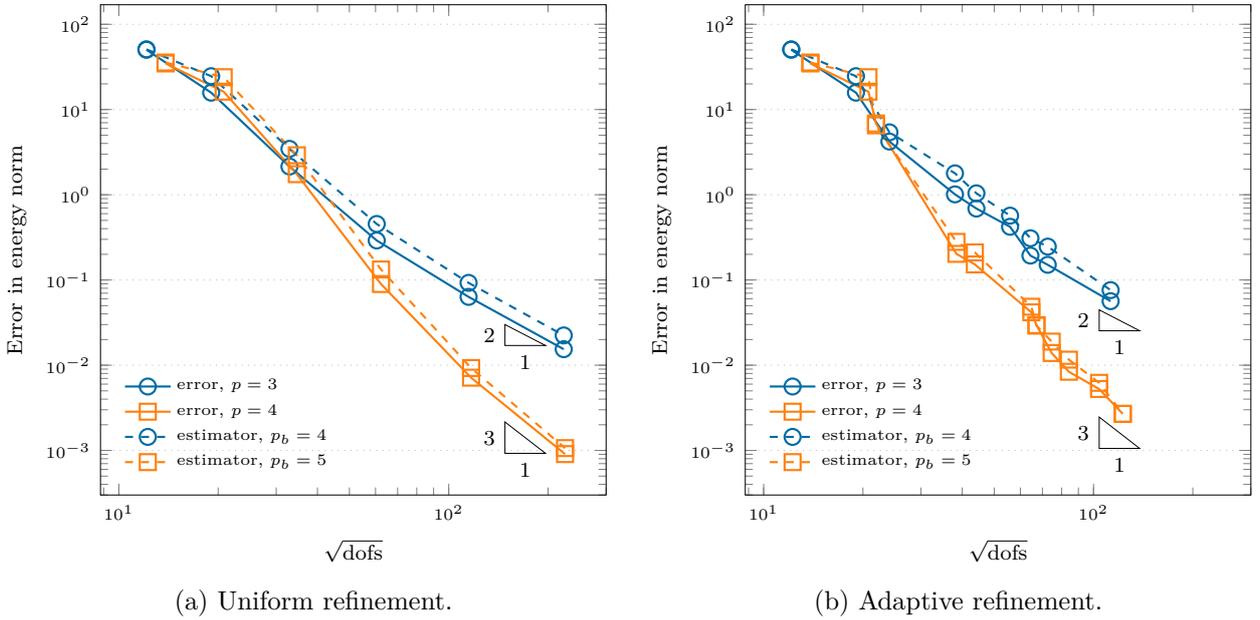

	\centering
	\begin{subfigure}[t]{0.495\textwidth}
	\input{\graphDir/testcase_KL_scordelis_with_hole_bubbles_GR.tex}
		\caption{Uniform refinement.}
	\end{subfigure}
	\hfill
	\begin{subfigure}[t]{0.495\textwidth}
	\input{\graphDir/testcase_KL_scordelis_with_hole_bubbles_MS.tex}
		\caption{Adaptive refinement.}
	\end{subfigure}	
	\caption{Study of the convergence of the error measured in the (approximated) energy norm and the bubble estimator for the Scordelis-Lo roof with a hole defined in test case 1. Comparison of uniform and adaptive refinements.}
\label{fig:convergence_scordelis_with_hole}
\end{figure}

In this numerical example, we take the geometry and material properties as defined in the well-known Scordelis-Lo roof benchmark, for a detailed description we refer to \citep{Belytschko1985}. We impose rigid diaphragm boundary conditions at both curved ends of the cylindrical structure, meaning that we hinder the displacement in the $xz$-plane as defined by the coordinate system of~\Cref{fig:scordelis_with_hole_geo}. Then, we trim out a circular hole in the parameter space of the surface, as depicted in~\Cref{fig:scordelis_with_hole_geo}. We recall that, analogously to the original setup, the structure is subjected to its self-weight and we set its value to $f_{z} = -90 \, [N/m^2]$. To the best of the authors' knowledge there is no analytical solution in closed form for this problem, therefore we perform a convergence study of the error measured in the energy norm against a reference solution obtained with bi-quintic B-splines defined on a uniform mesh with $62,456$ active elements.
We define this quantity as:
\begin{align*}
\vert \vert \tilde{e} \vert \vert_{E(\Omega)} = \sqrt{a(\boldsymbol{u_h^{\text{ref}}} - \boldsymbol{u_h},\boldsymbol{u_h^{\text{ref}}} - \boldsymbol{u_h})} \sim \vert \vert \boldsymbol{u} - \boldsymbol{u_h} \vert \vert_{E(\Omega)} \, ,
\end{align*} 
where $\boldsymbol{u_h^{\text{ref}}}$ denotes the reference solution.
The results are depicted in~\Cref{fig:convergence_scordelis_with_hole} for THB-splines of degree $p=3,4$, where uniform and adaptive refinement are analyzed. In both cases, we observe optimal rates of convergence of the error and a similar level of accuracy between the uniform and the adaptive strategy. This is due to the fact that this problem exhibits a very smooth solution. 

\subsubsection{Test case 2: Scordelis-Lo roof with holes under a point load}

\begin{figure}[!ht]
	\centering
	\includegraphics[width=0.85\textwidth]{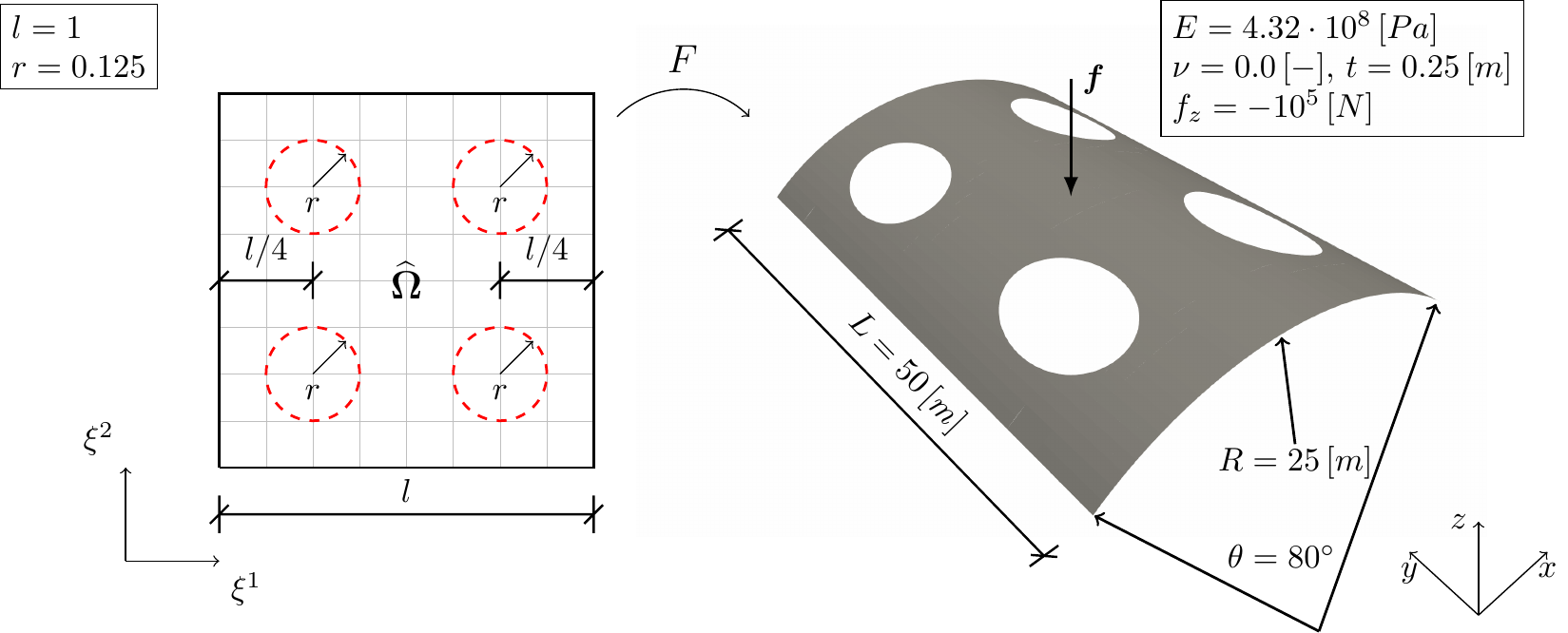}
	\caption{Geometry and physical parameters of the Scordelis-Lo roof with four holes.} \label{fig:scordelis_with_holes_geo}
\end{figure}

\begin{figure}
		\centering
		\begin{subfigure}[t]{0.325\textwidth}
			\centering
			\includegraphics[width=\textwidth]{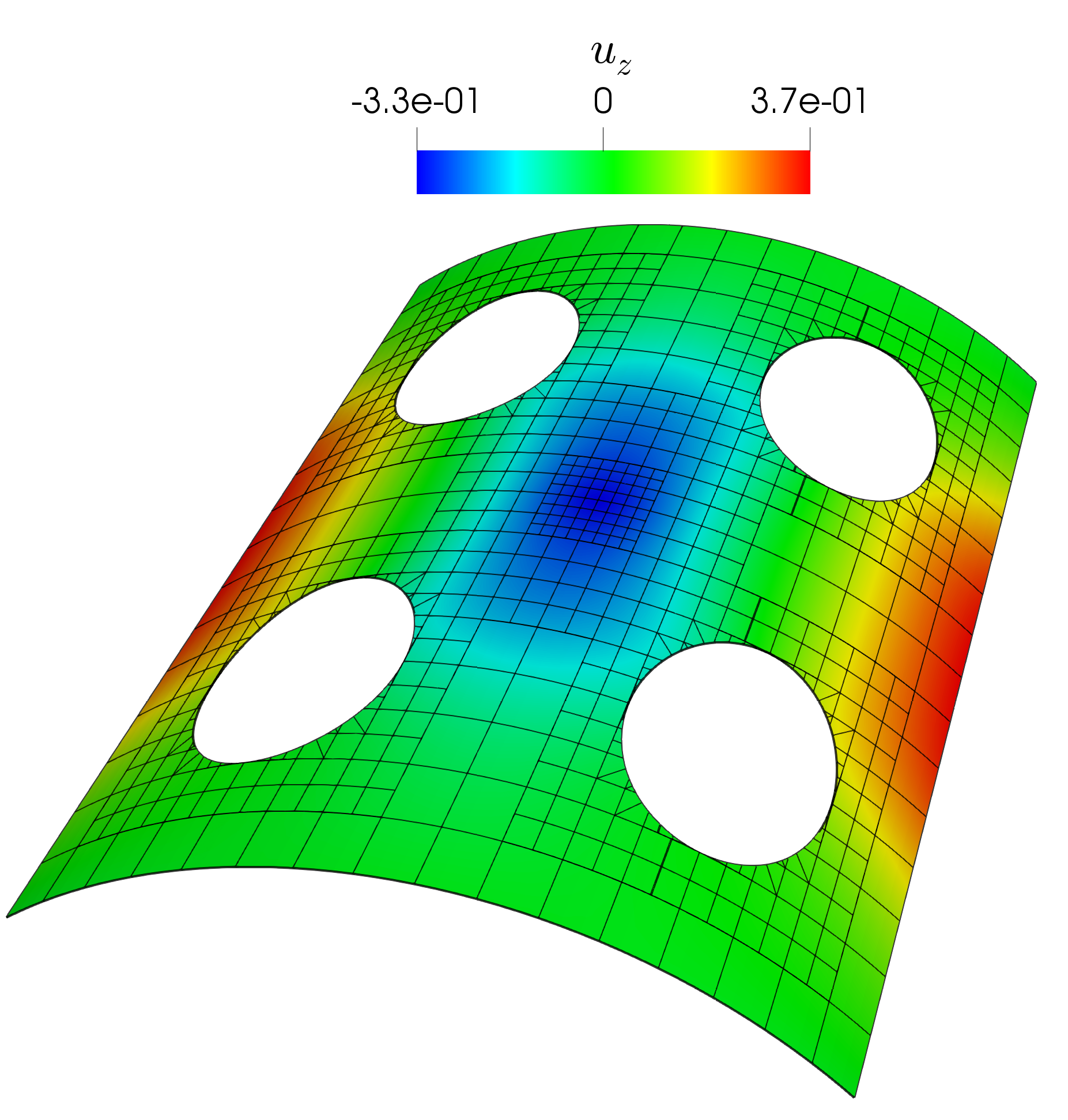}
			\caption{Mesh and $u_z$ at iteration $k = 3$.} \label{fig:scordelis_with_holes_a}
		\end{subfigure}
		\begin{subfigure}[t]{0.325\textwidth}
			\centering
			\includegraphics[width=\textwidth]{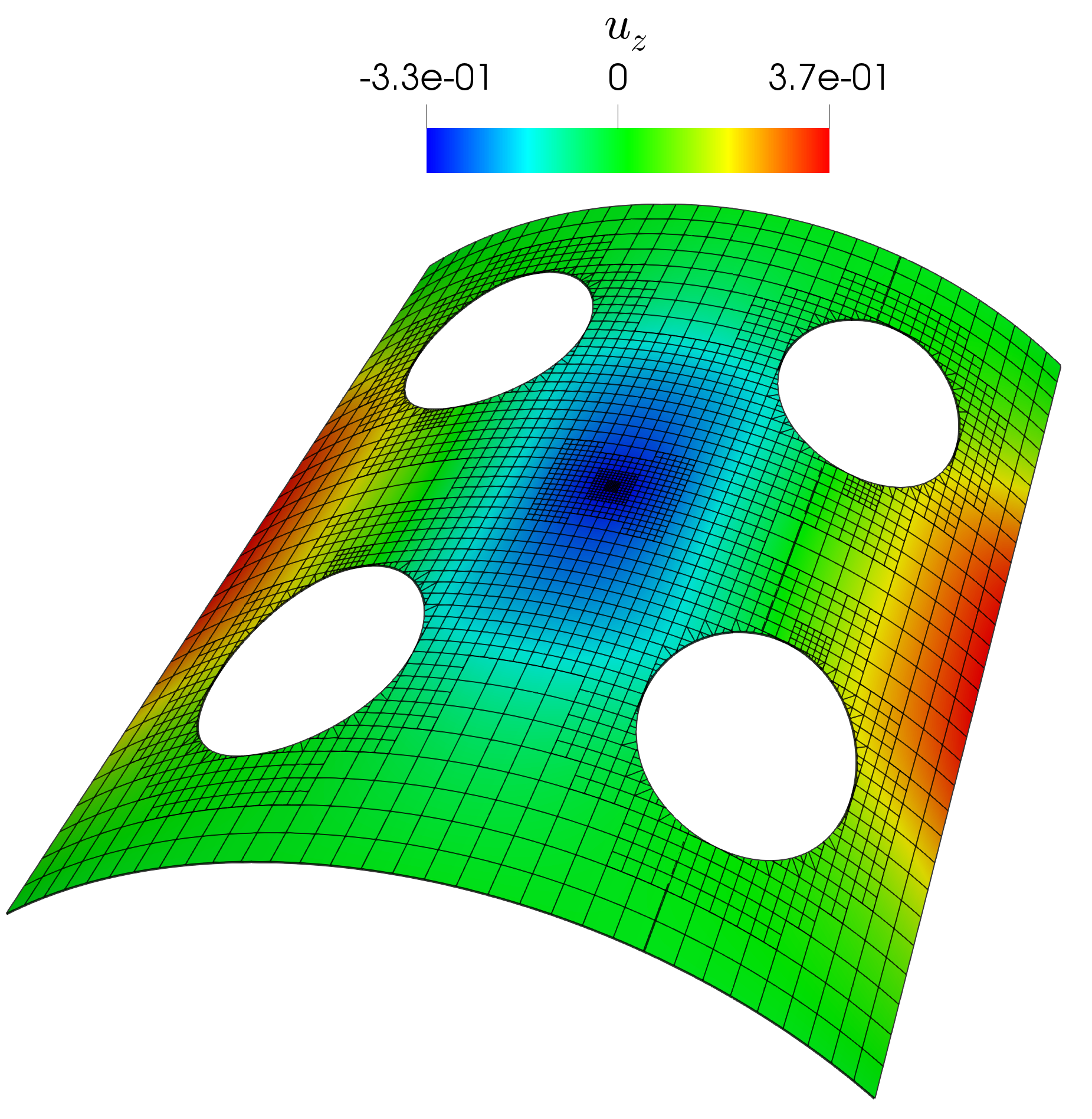}
			\caption{Mesh and $u_z$ at iteration $k = 6$.} \label{fig:scordelis_with_holes_b}
		\end{subfigure}
		\begin{subfigure}[t]{0.325\textwidth}
			\centering
			\includegraphics[width=\textwidth]{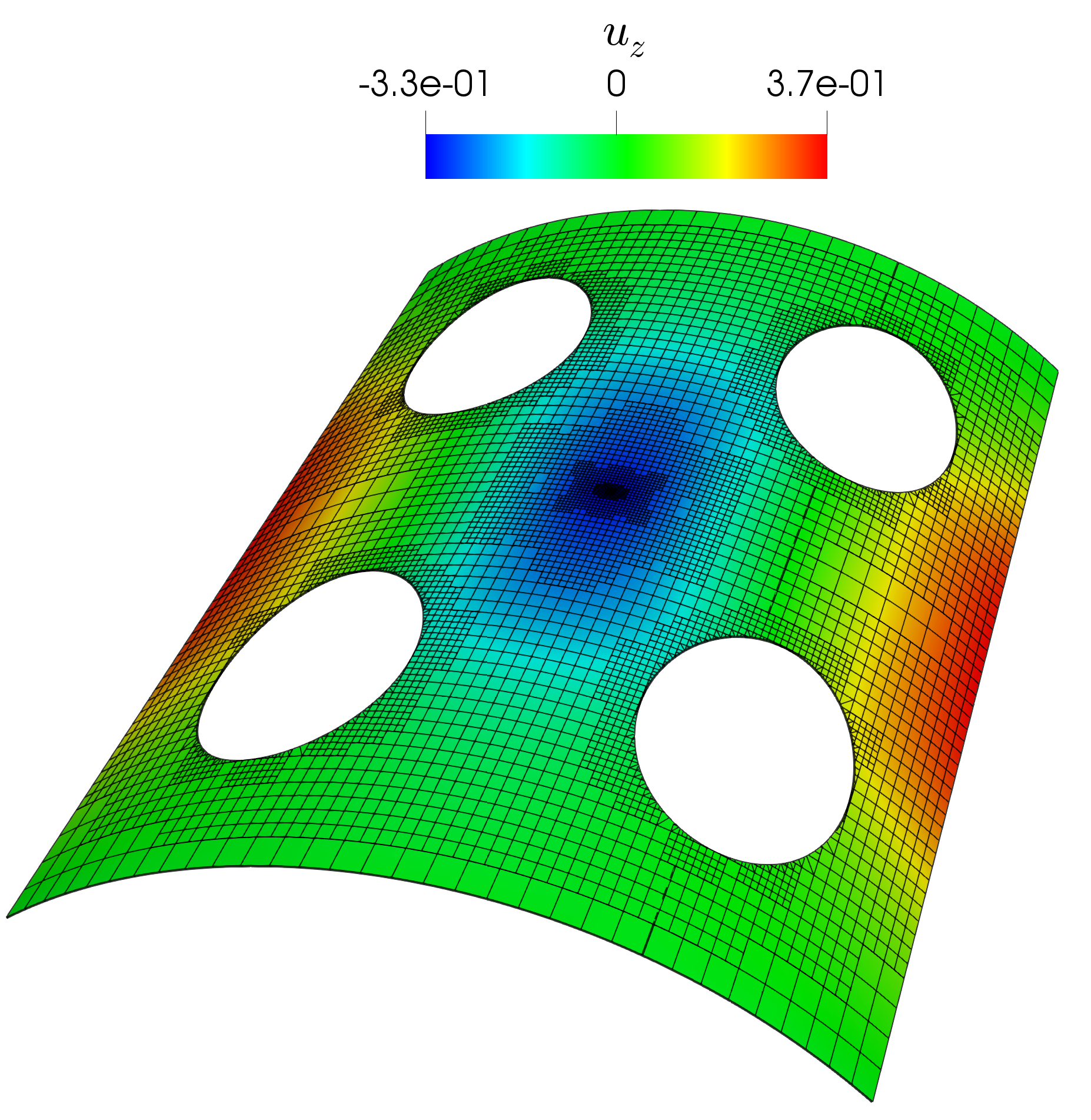}
			\caption{Mesh and $u_z$ at iteration $k = 8$.} \label{fig:scordelis_with_holes_c}
		\end{subfigure}
		\begin{subfigure}[t]{0.325\textwidth}
			\centering
			\includegraphics[width=\textwidth]{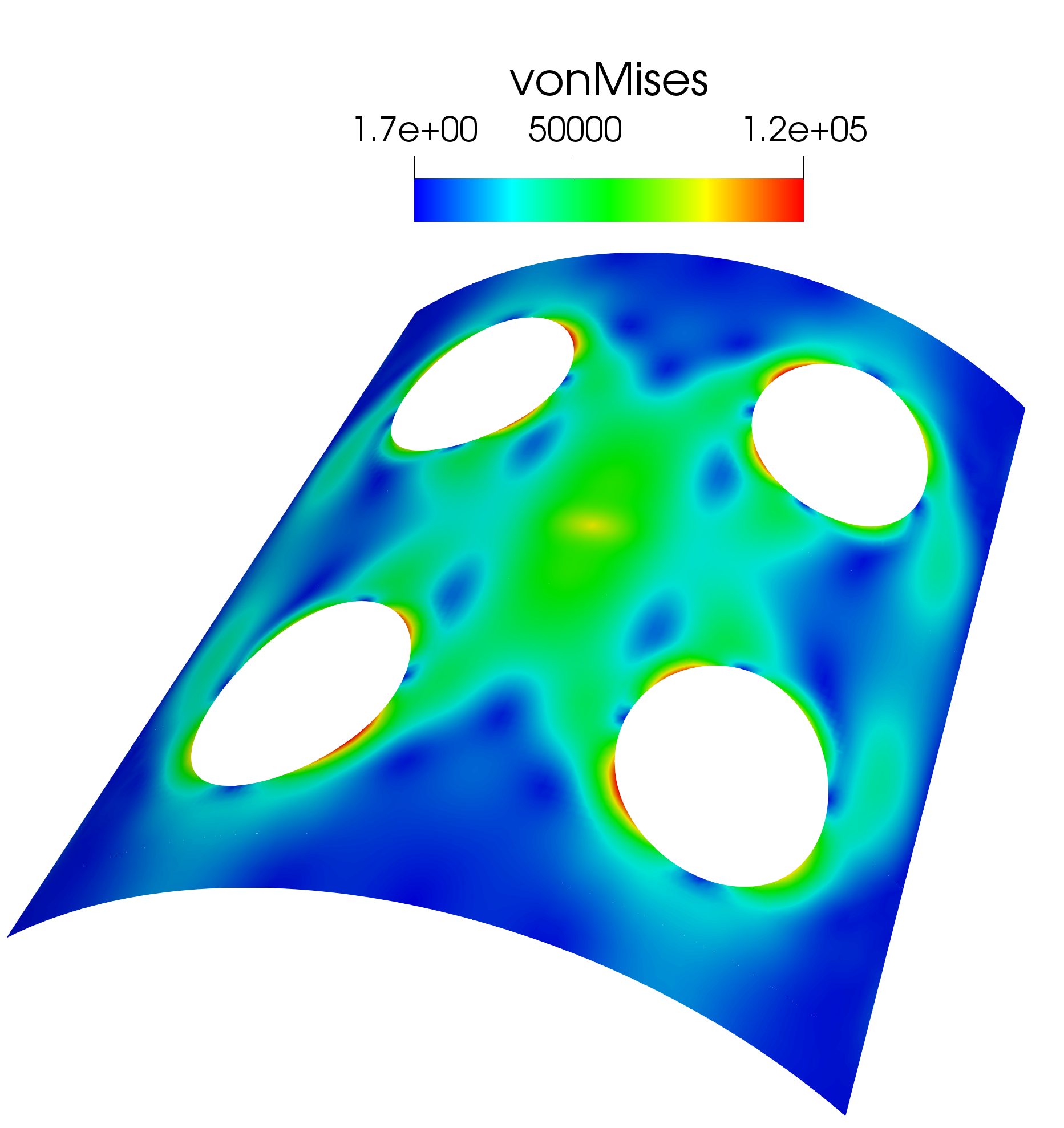}
			\caption{Von Mises at iteration $k = 3$.} \label{fig:scordelis_with_holes_d}
		\end{subfigure}
		\begin{subfigure}[t]{0.325\textwidth}
			\centering
			\includegraphics[width=\textwidth]{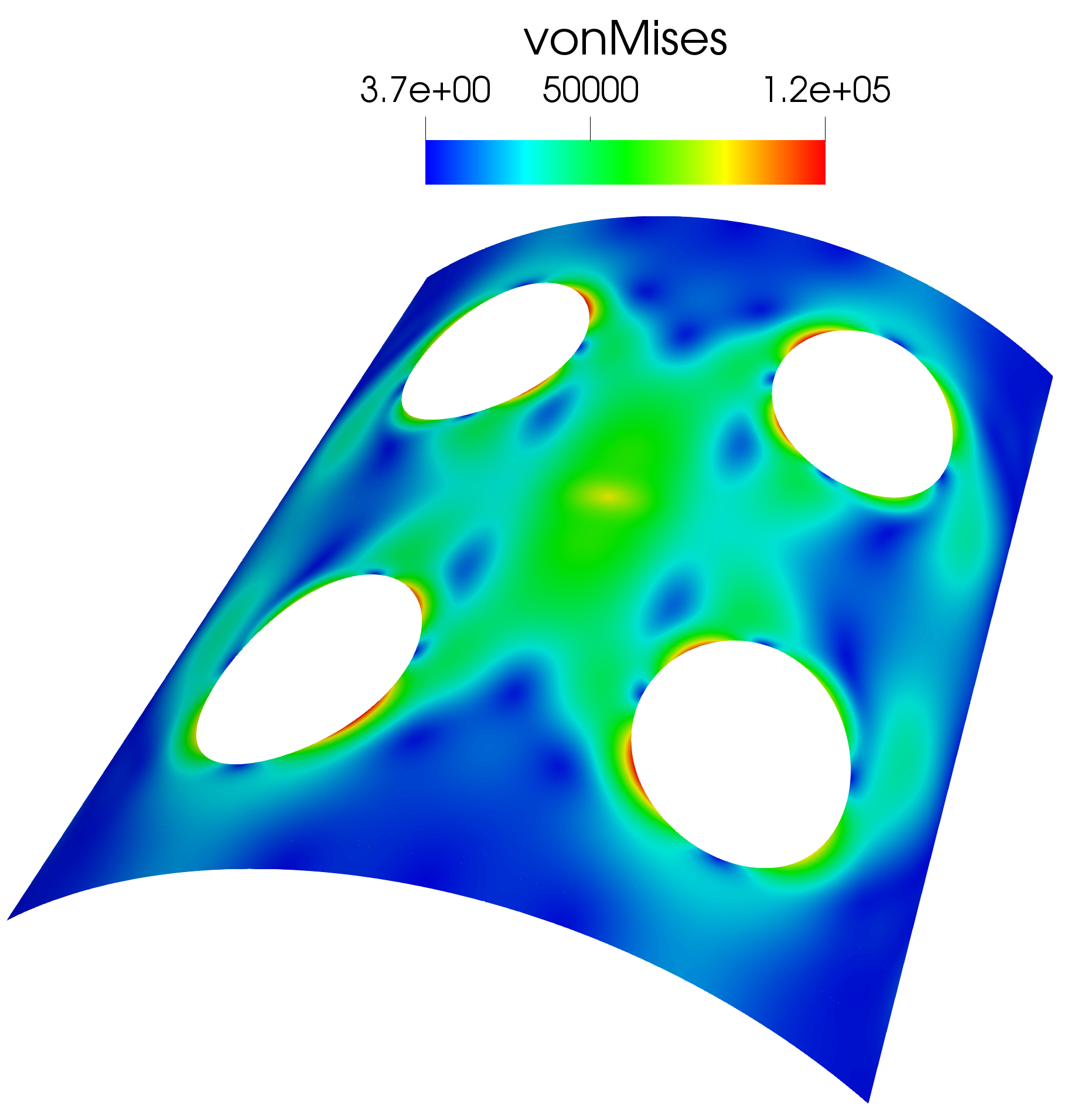}
			\caption{Von Mises at iteration $k = 6$.} \label{fig:scordelis_with_holes_e}
		\end{subfigure}
		\begin{subfigure}[t]{0.325\textwidth}
			\centering
			\includegraphics[width=\textwidth]{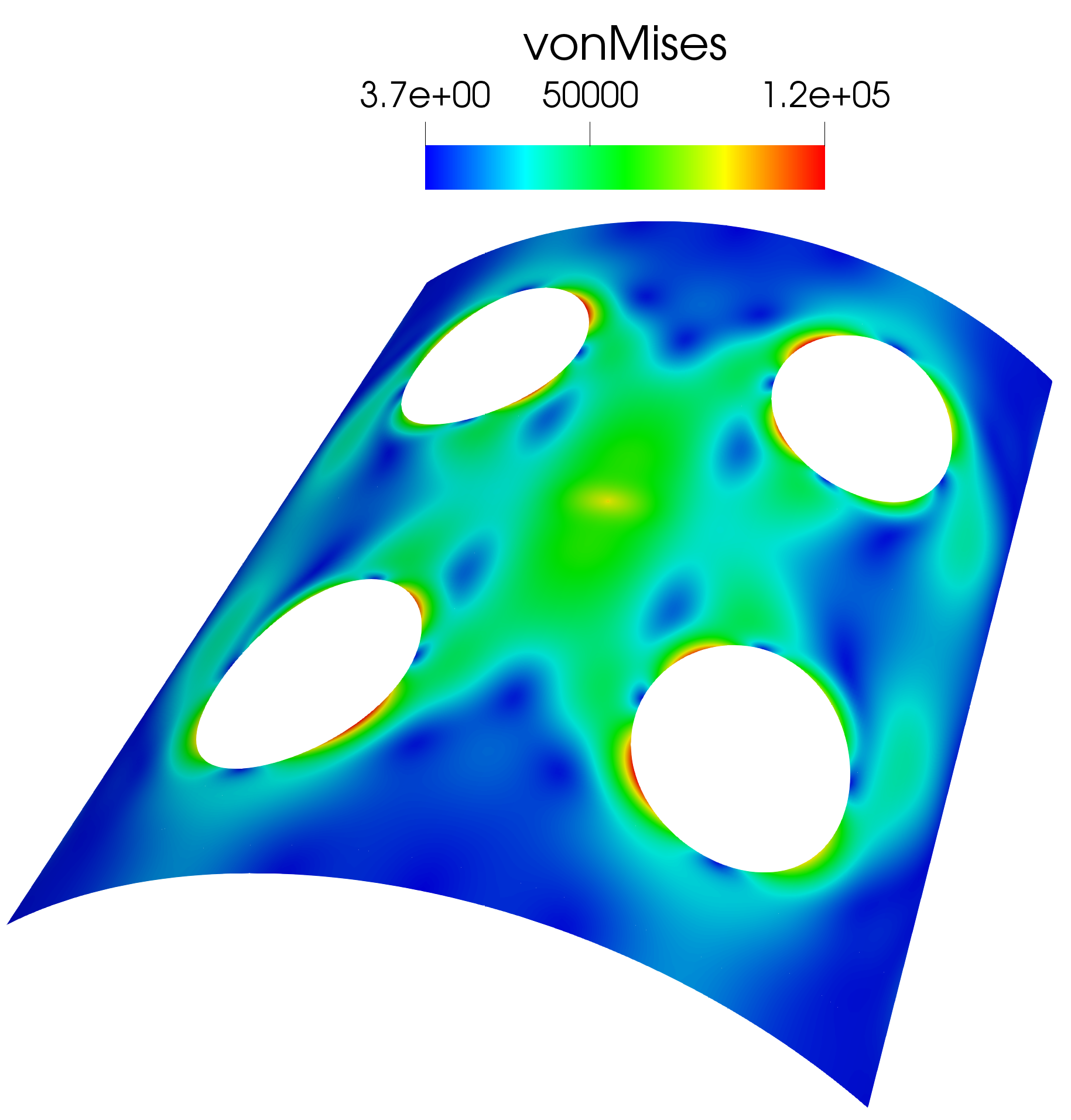}
			\caption{Von Mises at iteration $k = 8$.} \label{fig:scordelis_with_holes_f}
		\end{subfigure}		
		\par\bigskip
		\begin{subfigure}[t]{0.7\textwidth}
			\centering
			\input{\graphDir/convergence_KL_scordelis_with_holes_displacement.tex}
			\caption{Convergence of the normalized quantity $\vert 1 - u_z / u_z^{\text{ref}} \vert$.} \label{fig:scordelis_with_holes_convergence_displ}
		\end{subfigure}
					\caption{Mesh, solution and Von Mises stress at different steps of the adaptive loop driven by the bubble error estimator for the Scordelis-Lo roof with holes subjected to a point load, solution obtained employing truncated hierarchical B-splines of degree $p=3$ (a-f). Convergence plot of the normalized displacement at the center of the structure. Comparison of uniform and adaptive refinements (g).}
		\label{fig:scordelisWithHoles_mesh_and_vonMises}
\end{figure}

In this numerical test, we use again the same geometrical and material properties as defined in the Scordelis-Lo roof benchmark. This time, we trim the structure with four holes defined in the parametric space of the surface as circles, as depicted in~\Cref{fig:scordelis_with_holes_geo}. As in the example above, both curved ends are supported by rigid diaphragms boundary conditions (the displacement in the $xz$-plane is set to zero, with respect to the coordinate system defined in~\Cref{fig:scordelis_with_holes_geo}). Additionally, the roof is subjected to a point load in the vertical direction, defined as $f_z = - 10^5 \, [N]$, applied at the center of the structure. 
A reference displacement $u_z^{\text{ref}}$ at the location under the point load has been obtained with a fine, uniformly refined mesh of bi-cubic B-splines with $53,216$ active elements. In~\Cref{fig:scordelis_with_holes_convergence_displ} the convergence results for the uniform and the adaptive refinement are compared, where on the $y-$axis we plot the quantity $\vert 1 - u_z / u_z^{\text{ref}} \vert$, evaluated at the center of the structure. We notice that the adaptive strategy achieves a solution that is several order of magnitudes more accurate than uniform refinement for the same number of degrees of freedom. This result confirms the superior efficiency of our error-driven strategy in the presence of singularities. Lastly, in~\Cref{fig:scordelis_with_holes_a,fig:scordelis_with_holes_b,fig:scordelis_with_holes_c,fig:scordelis_with_holes_d,fig:scordelis_with_holes_e,fig:scordelis_with_holes_f} we present the numerical solution $u_z$ and corresponding Von Mises stress at different iterations $k=3,6,8$ of the adaptive simulation. Remarkably, the refinement is driven not only in the region where the load is applied but also in those areas around the holes where stress concentrations are present.   	 

%
%

\subsubsection{Test case 3: B-pillar of a car}

\begin{figure}
		\centering
		\begin{subfigure}[t]{0.35\textwidth}
			\centering
			\includegraphics[width=\textwidth]{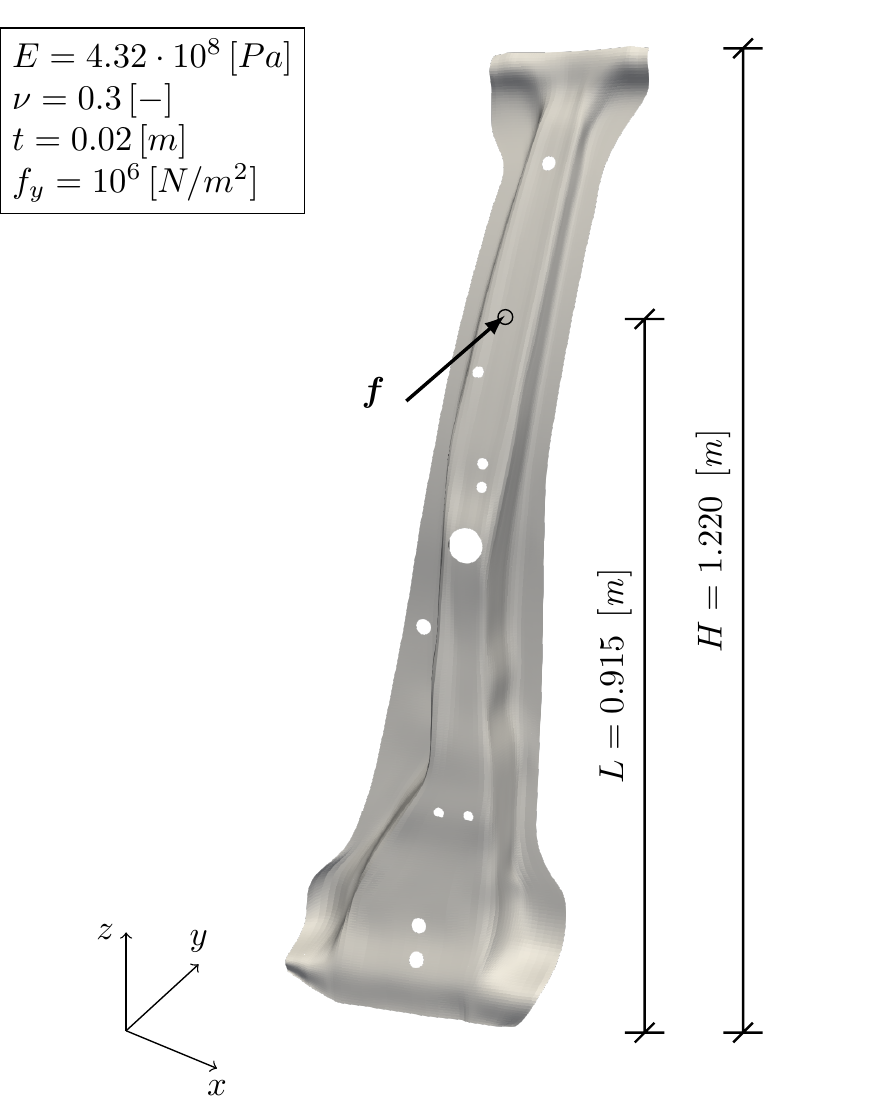}
			\caption{Problem setup.}\label{fig:b_pillar_geo}
		\end{subfigure}		
		\begin{subfigure}[t]{0.30\textwidth}
			\centering
			\includegraphics[width=\textwidth]{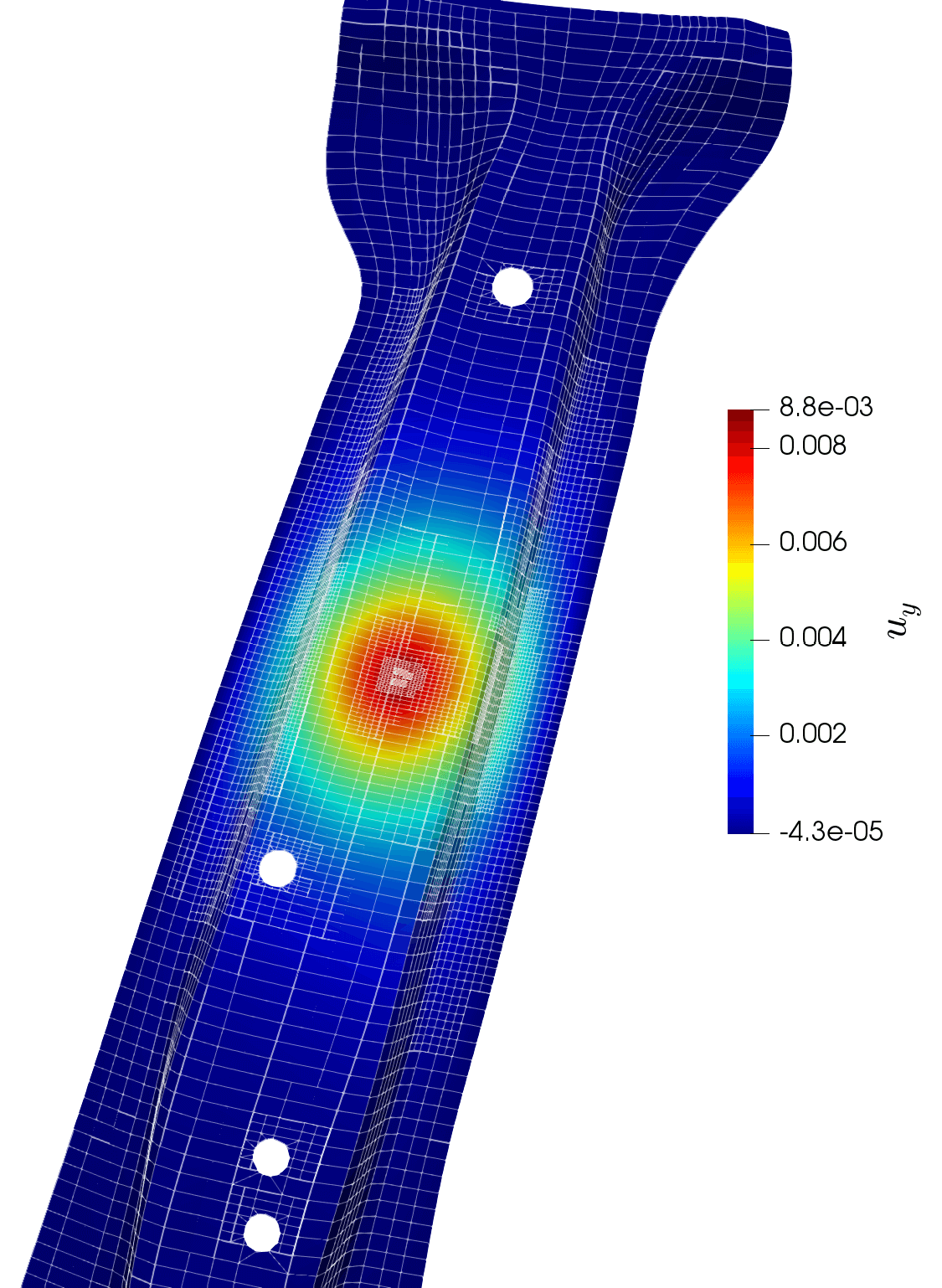}
			\caption{Mesh and solution $u_y$ at iteration $k = 5$.}\label{fig:b_pillar_mesh_and_vonMises1}
		\end{subfigure}
		\begin{subfigure}[t]{0.30\textwidth}
			\centering
			\includegraphics[width=\textwidth]{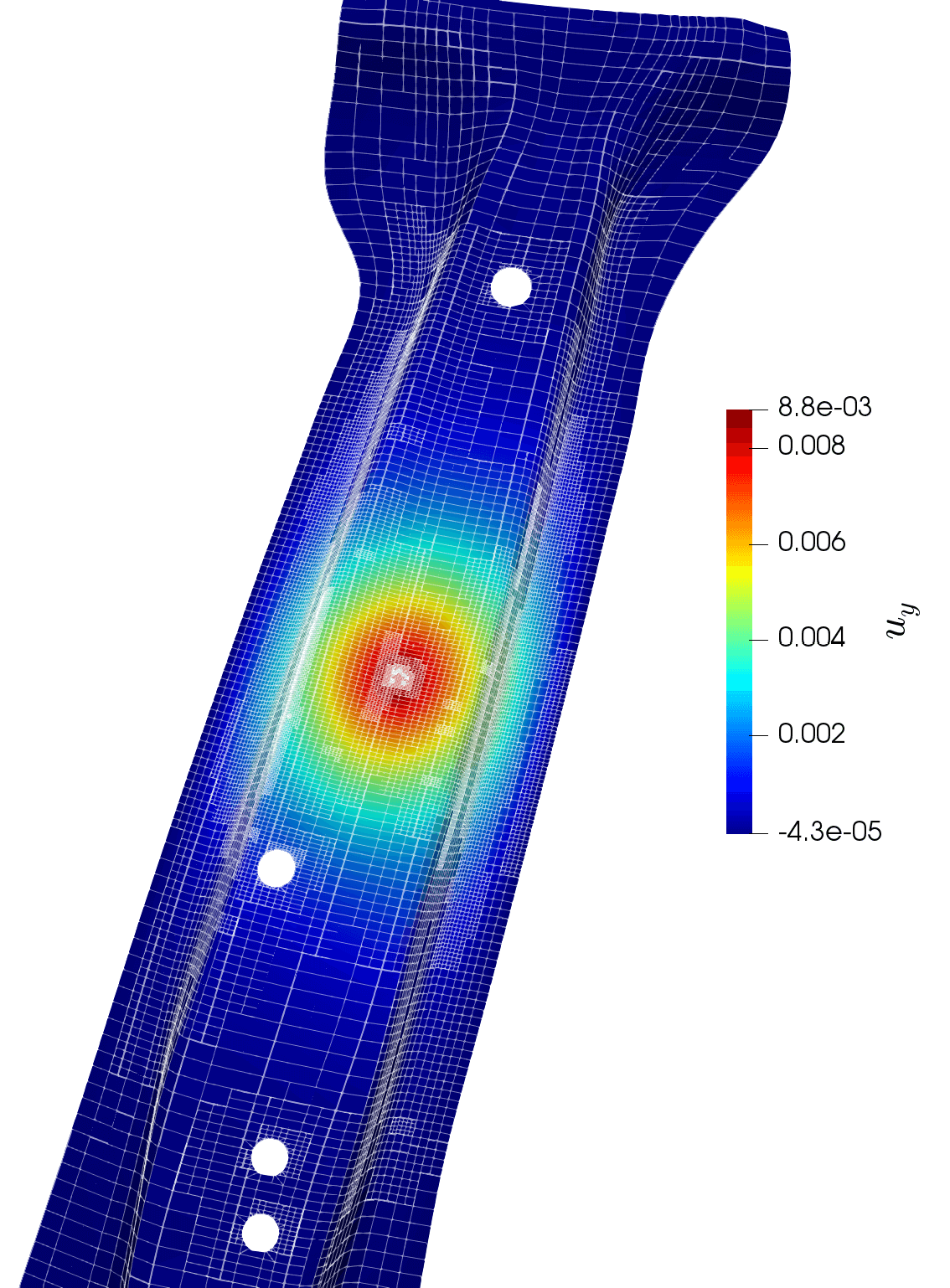}
			\caption{Mesh and solution $u_y$ at iteration $k = 7$.}\label{fig:b_pillar_mesh_and_vonMises2}
		\end{subfigure}
		\begin{subfigure}[t]{0.30\textwidth}
			\centering
			\includegraphics[width=\textwidth]{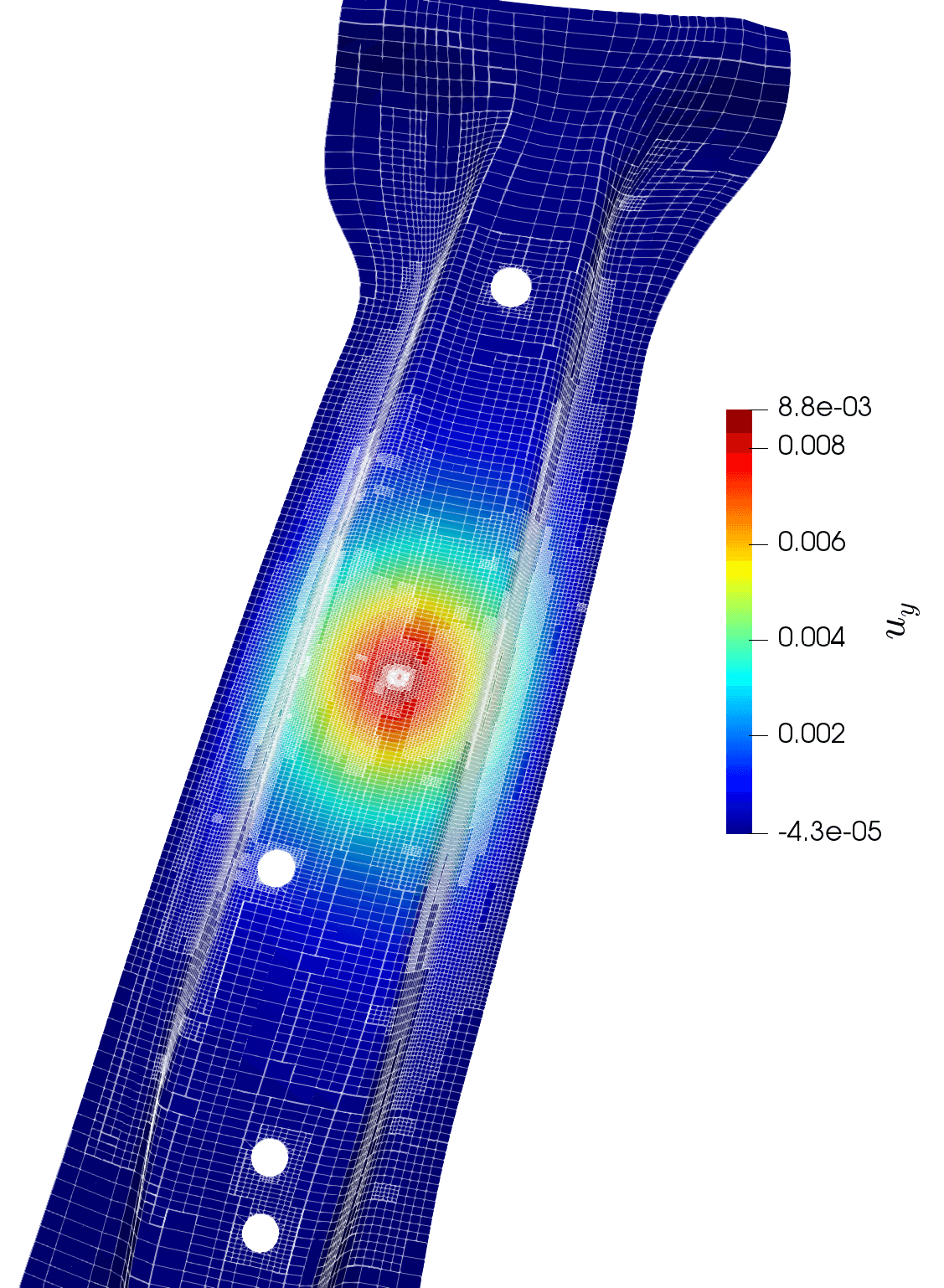}
			\caption{Mesh and solution $u_y$ at iteration $k = 9$.}\label{fig:b_pillar_mesh_and_vonMises3}
		\end{subfigure}
		\begin{subfigure}[t]{0.385\textwidth}
			\centering
			\includegraphics[width=\textwidth]{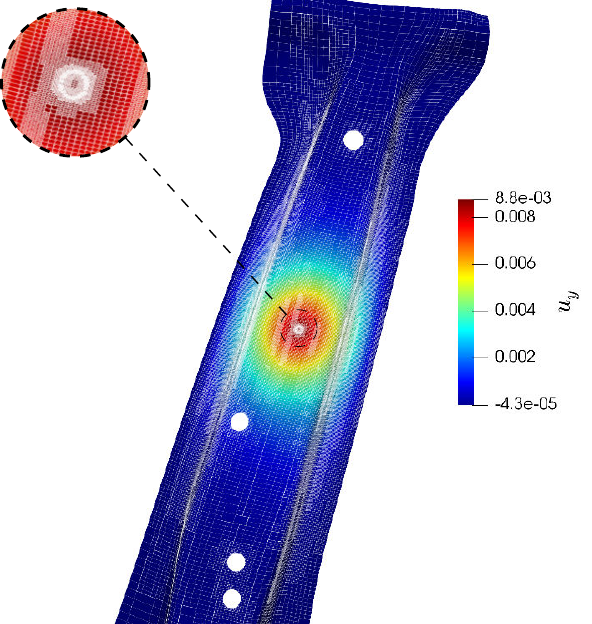}
			\caption{Mesh and solution $u_y$ at iteration $k = 11$.}\label{fig:b_pillar_mesh_and_vonMises4}
		\end{subfigure}
		\begin{subfigure}[t]{0.295\textwidth}
			\centering
			\includegraphics[width=\textwidth]{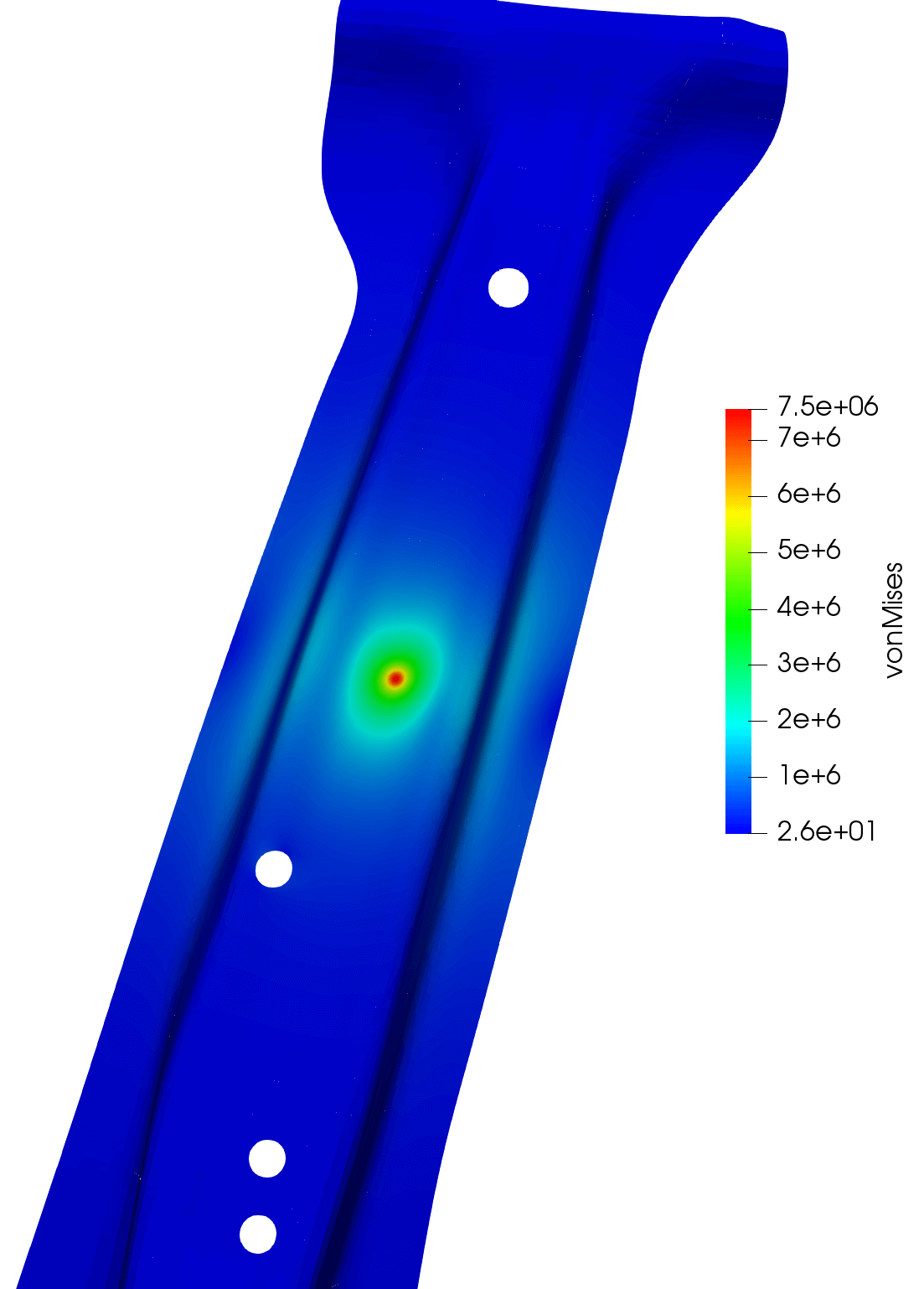}
			\caption{Von Mises at iteration $k = 11$.}\label{fig:b_pillar_mesh_and_vonMises5}
		\end{subfigure}
		\caption{Geometry of the B-pillar and mesh, solution and Von Mises stress at different steps of the adaptive loop driven by the bubble error estimator for the B-pillar example, solution obtained employing truncated hierarchical B-splines of degree $p=3$.}
		\label{fig:b_pillar_mesh_and_vonMises}
\end{figure}

In the last example, we model the B-pillar of a car in the commercial CAD software \textit{Rhinoceros}\footnote{http://www.rhino3d.com}. The geometry is defined as a trimmed, single-patch B-spline surface of degree $p=3$ composed of 30 and 180 knot spans in the two parametric directions, respectively. The structure is considered fixed on the entire external boundary (all displacement components are set to zero) and is loaded with a circular, uniformly distributed force of radius $r = 0.01 \, [m]$ in the y-direction, see~\Cref{fig:b_pillar_geo}.
This example demonstrates the applicability of our approach to the adaptive analysis of complex shell structures of industrial relevance, modeled and exported directly from a commercial CAD software. In~\Cref{fig:b_pillar_mesh_and_vonMises1,fig:b_pillar_mesh_and_vonMises2,fig:b_pillar_mesh_and_vonMises3,fig:b_pillar_mesh_and_vonMises4,fig:b_pillar_mesh_and_vonMises5} the numerical solution at several steps $k=5,7,9,11$ of the adaptive loop and the Von Mises stress obtained with truncated hierarchical B-splines of degree $p=3$ are depicted. It is worth noting how the estimator correctly refines the circular area around the load, particularly at the boundary of the circle where a sharp change in boundary conditions is present. Furthermore, it also detects and refines stress concentrations in the proximity of the trimmed holes and regions of high curvature of the geometry, where higher bending stresses are expected.

\section{Conclusions} \label{sec:conclusions}

We have presented an innovative adaptive framework in the scope of second- and fourth-order elliptic partial differential equations defined on trimmed domains, which makes use of the re-parametrization tool for the integration of trimmed elements introduced in~\citep{Antolin2019} and exploits the local refinement capabilities of truncated hierarchical B-splines. Related to the latter, we highlight that truncated hierarchical B-splines retain their linear independence also when the domain on which they are defined is cut, where a simple algorithm is given which preserves this property in standard isogeometric implementations. 
For steering the automatic refinement of the mesh, we introduced an a posteriori error estimator which was inspired by some previous works in the field of multi-level estimators \citep{Bank1993,Vuong2011} and it is a direct extension of the method presented by the authors in~\citep{Antolin2019b}.
We recall that the estimator is based on the solution of an additional, residual-like system. As presented in details in~\citep{Antolin2019b}, the resulting linear system has in general the following properties:
\begin{enumerate}

\item block-diagonal structure due to the single element support of bubble functions,

\item small dimension of each block, that are defined element-wise,


\item easily parallelizable.

\end{enumerate}
All the aforementioned properties make the proposed method appealing from a computational point of view and easy to implement into existing isogeometric codes. 

To the best of the authors' knowledge, this is the first methodology that allows an error-driven simulation in the context of trimmed isogeometric Kirchhoff-Love shells. In particular, it is a step towards the integration of design and automated adaptive isogeometric analysis, where, at first, we have validated the proposed framework on an extensive series of benchmark problems, where we have always observed optimal rates of convergence, and later we have demonstrated its applicability to complex, industrial-like geometries. This was achieved thanks to the interoperability between trimmed surfaces created in the commercial CAD software 
\textit{Rhinoceros}
and the in-house tool for handling cut elements~\citep{Antolin2019}.     

To conclude, we have numerically proved the reliability and robustness of the proposed error-driven adaptive framework for the simulation of various mechanically-relevant PDEs defined on trimmed domains, namely in the scope of the Poisson problem, linear elasticity and Kirchhoff-Love shells. We have shown that our method systematically achieves superior efficiency and accuracy per-degree-of-freedom in problems where sharp features of the solution and/or singularities are present, where we exploit truncated hierarchical B-splines to achieve local refinement of the basis.
Finally, we have successfully applied it to the Kirchhoff-Love shell analysis of a complex trimmed surface created in \textit{Rhinoceros}, namely the B-pillar of a car, which demonstrate the potential of our framework to be readily used to tackle problems of engineering relevance.

\section*{Acknowledgements} 
The authors would like to thank Prof. Alessandro Reali for the fruitful discussion on the subject of this paper and Kenji Takada from \textit{Honda R}\&\textit{D Co.}\footnote{https://global.honda} for providing the original CAD geometry of the B-pillar. Furthermore, we would like to thank Dr. Xiaodong Wei for his help with the pre-processing of the geometry.
The authors also gratefully acknowledge the support of the European Research Council, via the ERC AdG project CHANGE n.694515.

\bibliographystyle{plainnatnourl}
\bibliography{library.bib}

\end{document}